\documentclass{amsart}
\usepackage{amssymb} 
\usepackage{amsmath} 
\usepackage{amscd}
\usepackage{amsbsy}
\usepackage{comment}
\usepackage[matrix,arrow]{xy}
\usepackage{hyperref}

\DeclareMathOperator{\Norm}{Norm}

\DeclareMathOperator{\Ker}{Ker}

\DeclareMathOperator{\Res}{Res}

\DeclareMathOperator{\Gal}{Gal}
\DeclareMathOperator{\norm}{Norm}

\DeclareMathOperator{\ord}{ord}

\DeclareMathOperator{\lcm}{lcm}

\newcommand{\Q}{{\mathbb Q}}
\newcommand{\Z}{{\mathbb Z}}

\newcommand{\R}{{\mathbb R}}
\newcommand{\F}{{\mathbb F}}

\newcommand{\cM}{\mathcal{M}}
\newcommand{\cY}{\mathcal{Y}}
\newcommand{\cE}{\mathcal{E}}
\newcommand{\cF}{\mathcal{F}}
\newcommand{\cN}{\mathcal{N}}

\newcommand{\cQ}{\mathcal{Q}}
\newcommand{\cR}{\mathcal{R}}
\newcommand{\cT}{\mathcal{T}}
\newcommand{\cS}{\mathcal{S}}

\newcommand{\OO}{{\mathcal O}}

\newcommand{\ff}{\mathfrak{f}}
\newcommand{\fg}{\mathfrak{g}}
\newcommand{\fp}{\mathfrak{p}}
\newcommand{\fq}{\mathfrak{q}}

\newcommand{\bs}{\mathbf{s}}

\def\mod#1{{\ifmmode\text{\rm\ (mod~$#1$)}
\else\discretionary{}{}{\hbox{ }}\rm(mod~$#1$)\fi}}

\begin {document}

\newtheorem{thm}{Theorem}
\newtheorem{lem}{Lemma}[section]
\newtheorem{prop}[lem]{Proposition}

\theoremstyle{definition}

\theoremstyle{remark}

\title[]{Shifted powers in Lucas-Lehmer sequences}
%\author{Samir Siksek}
%\address{Mathematics Institute\\
%	University of Warwick\\
%	Coventry\\
%	CV4 7AL \\
%	United Kingdom}

%\email{s.siksek@warwick.ac.uk}
\author[Michael Bennett]{Michael A. Bennett}
\address{Department of Mathematics, University of British Columbia, Vancouver, B.C., V6T 1Z2 Canada}
\email{bennett@math.ubc.ca}

\author{Vandita Patel}
\address{Department of Mathematics, University of Toronto, Toronto, Canada, M5S 2E4}
\email{vandita@math.utoronto.ca}

\author{Samir Siksek}
\address{Mathematics Institute, University of Warwick, Coventry CV4 7AL, United Kingdom}
\email{S.Siksek@warwick.ac.uk}
\thanks{ The first-named is supported by NSERC. The third-named author is supported by an EPSRC Leadership Fellowship EP/G007268/1, and EPSRC {\em LMF: L-Functions and Modular Forms} Programme Grant EP/K034383/1.  }

\date{\today}

\keywords{Exponential equation,
Lucas sequence, shifted power, Galois representation,
Frey curve,
modularity, level lowering, Baker's bounds, Hilbert modular forms,
Thue equation}
\subjclass[2010]{Primary 11D61, Secondary 11D41, 11F80, 11F41}

\begin {abstract}
We develop a general framework for finding all perfect powers in sequences derived by shifting non-degenerate quadratic Lucas-Lehmer binary recurrence sequences by a fixed integer. By combining this setup with bounds for linear forms in logarithms and results based upon the modularity of elliptic curves defined over totally real fields, we are able to answer a question of Bugeaud, Luca, Mignotte and the third author by explicitly finding all perfect powers of the shape $F_k \pm 2 $ where $F_k$ is the $k$-th term in the Fibonacci sequence.
\end {abstract}
\maketitle

%------------------------------
\section{Introduction}
%------------------------------

If $\{ u_n \}$ is a non-degenerate integer binary linear recurrence sequence, then the sequence $\{ u_n \}$ contains at most finitely many integer perfect powers, which may be effectively determined. This result was 	proved independently, using bounds for linear forms in Archimedean and non-Archimedean logarithms, by Peth\H{o} \cite{Pe} and Shorey and Stewart \cite{ShSt1}. The explicit determination of all such powers in a given sequence, however, has been achieved in only a few cases, principally in those where the problem may be reduced to a question of solving ternary Diophantine equations with integer coefficients. In such a situation, the possibility exists to combine the machinery of linear forms in logarithms with information derived from considering certain Frey--Hellegouarch curves corresponding to the ternary equations. A prototype for these problems may be found in the paper of Bugeaud, Mignotte and the third author \cite{BMS}, where all perfect powers in the Fibonacci sequence are determined; this amounts to finding the integer solutions to the equation
$$
x^2 - 5 y^{2p} = \pm 4,
$$
in prime numbers $p$ and integers $x$ and $y$. Here, results from the theory of linear forms in logarithms provide a manageable upper bound upon the exponent $p$, but solving the remaining (hyperelliptic) equations is accomplished only through considering them as ternary equations of signature $(p,p,2)$ and using arguments based upon the modularity of Galois representations to deduce arithmetic information guaranteeing that $x$ is necessarily extraordinarily large (unless $x = \pm 1$).

If we shift a given recurrence, considering, say, $u_n+c$ for a nonzero integer $c$, instead of just $u_n$, the situation becomes considerably  more complicated. The resulting sequence need not possess much of the basic structure of a binary linear recurrence sequence, despite sharing a similar rate of growth. In particular, various divisibility statements may no longer hold, and questions of the existence of primitive divisors are significantly harder to address. Despite this, Shorey and Stewart \cite{SS} were able to show, under mild hypotheses,
 that, given fixed integers $a$ and $c$, the equation
$$
u_n+c = a y^p
$$ 
has at most finitely many, effectively computable solutions. Only in very special cases, however, can such equations correspond to Frey--Hellegouarch curves defined over $\mathbb{Q}$ (see e.g.  the paper of Bugeaud, Luca, Mignotte and the third author \cite{BLMS} for a number of such examples).

In a previous paper \cite{BDMS}, the first and third authors, with Dahmen and Mignotte, developed a method combining information derived from Frey--Hellegouarch curves defined over real quadratic fields with lower bounds for linear forms in logarithms to explicitly determine all shifted powers in certain binary recurrence sequences. The setup in \cite{BDMS} was as follows. 
Let $K$ be a real quadratic number field, $\mathcal{O}_K$ its ring of integers and $\varepsilon \in \mathcal{O}_K$ a fundamental unit in $K$, with conjugate $\overline{\varepsilon}$. Define the Lucas sequences $U_k$ and $V_k$, of the first and second kinds, respectively, via
$$
U_k=\frac{\varepsilon^k-\left(\overline{\varepsilon}\right)^k}{\varepsilon-\overline{\varepsilon}} \; \mbox{ and } \;
V_k =  \varepsilon^k+\left(\overline{\varepsilon}\right)^k, \; \mbox{ for } \;  k \in \Z.
$$
Let $a,c \in \Z$ with $a\not=0$, and consider the problem of determining the shifted powers $ay^p-c$ in one of these sequences, i.e. determining all integers $k, y$ and $p$ with $p \geq 2$ prime (say) such that we have
\begin{equation}\label{eqn:ShiftedPowers}
U_k+c=ay^p
\end{equation}
or
\begin{equation}\label{eqn:ShiftedPowers2}
V_k+c=ay^p.
\end{equation}

In  \cite{BDMS},  techniques were introduced to potentially resolve such problems corresponding to either

\begin{itemize}
\item equation (\ref{eqn:ShiftedPowers}) with $k$ odd  and $\Norm(\varepsilon)=-1$, or
\item equation (\ref{eqn:ShiftedPowers2}) with either $k$ even or  $\Norm(\varepsilon)=1$.
\end{itemize}
Let us now describe an approach to treat the remaining cases.
For instance, a solution to  (\ref{eqn:ShiftedPowers}) leads to the equation
$$
\frac{\varepsilon^k-\left(\overline{\varepsilon}\right)^k}{\varepsilon-\overline{\varepsilon}} = a y^p-c
$$
and so we have
$$
\varepsilon^{2k} + (\varepsilon - \overline{\varepsilon} ) c \varepsilon^k-\Norm(\varepsilon)^k = 
(\varepsilon - \overline{\varepsilon} )  a \varepsilon^k y^p.
$$
It follows that
\begin{equation} \label{pp2-1}
\left( 2 \varepsilon^k + (\varepsilon - \overline{\varepsilon} )  c \right)^2 -
\left( 4 \Norm(\varepsilon)^k + (\varepsilon - \overline{\varepsilon} )^2 c^2 \right) =
4 (\varepsilon - \overline{\varepsilon} )  a \varepsilon^k y^p.
\end{equation}
Similarly, in the case of equation (\ref{eqn:ShiftedPowers2}), we have
$$ 
\varepsilon^k+\left(\overline{\varepsilon}\right)^k = a y^p-c
$$
whereby
\begin{equation} \label{pp2-2}
\left( 2 \varepsilon^k + c \right)^2 + 4 \Norm(\varepsilon)^k - c^2 = 4 a \varepsilon^k y^p.
\end{equation}
In either case, we can attach to a solution a Frey--Hellegouarch curve of signature $(p,p,2)$, defined over the totally real (quadratic) field $K$.

%---------------------------------------------------------------------
\section{Shifted powers in the Fibonacci sequence}
%--------------------------------------------------------------------

We will now describe an open question from the literature which our new techniques enable us answer.
Let $F_k$ be the {\it Fibonacci sequence} defined by
$$
F_0=0, \; \; F_1=1 \; \; \mbox{ and } \; \; F_{k+2}=F_{k+1}+F_k.
$$
Define further the {\it Lucas sequence} by
$$
L_0=2, \; \; L_1=1 \; \; \mbox{ and } \; \; L_{k+2}=L_{k+1}+L_k.
$$
For $K=\Q(\sqrt{5})$,  writing
$$
\varepsilon=\frac{1+\sqrt{5}}{2} \; \; \mbox{ and } \; \;
\overline{\varepsilon}=\frac{1-\sqrt{5}}{2},
$$
it follows that $\varepsilon$ is a fundamental unit of $K$ and, 
by Binet's formula,
$$
F_k=\frac{\varepsilon^k-\overline{\varepsilon}^k}{\sqrt{5}} 
\; \mbox{ and } \; 
L_k=\varepsilon^k+\overline{\varepsilon}^k,
$$
from which we obtain the well-known identity
\begin{equation}\label{eqn:wk}
L_k^2-5 F_k^2= 4 (-1)^k.
\end{equation}

In general, one has, for any integers $a$ and $b$,
\begin{equation} \label{stuff}
F_aL_b = F_{a+b} + (-1)^b F_{a-b}.
\end{equation}
This identity is used with $|a-b| \in \{ 1,2 \}$ in \cite{BLMS} to solve the
equations $F_k \pm 1 = y^p$ by reducing them to equations of the shape $F_k=
\alpha y^p$, for fixed integers $\alpha$ (which may be treated by considering
Frey--Hellegouarch curves defined over $\mathbb{Q}$). In this initial reduction,
it is of importance that $F_{-1}=F_1=F_2=1$ and $F_{-2}=-1$; more generally, 
%\margnote{As far as I can tell, the trick
%works for $F_n+c=y^p$ when $c=F_k$ where $k \equiv n \pmod{4}$. In this case we take $a=(n+k)/2$ and $b=(n-k)/2$.
%The congruence assumption guarantees that $b$ is even.}
analogous arguments allow one to treat equations of the form $F_n +c = y^p$,
for $c=F_{k}$ where $k \equiv n \mod{4}$. In particular, such a reduction does
not appear to be possible in general for the similar equation $F_k \pm 2 =y^p$
(which is posed an an open problem in \cite{BLMS}).

In this paper, we prove the following.
\begin{thm} \label{thm-main}
If $k, y$ and $p$ are integers, with $p$ prime and 
$$
F_k \pm 2 = y^p,
$$
then $|k| \in \{ 1, 2, 3, 4, 9\}$.
\end{thm}

Let us suppose that $F_k \pm 2=y^p$. In case $k$ is odd, say $k=2n+1$, choosing $a=n+2$ and $b=n-1$ in (\ref{stuff}), 
$$
F_k + (-1)^{n-1} 2 = F_{n+2}L_{n-1},
$$
while $a=n+1$ and $b=n-2$  gives
$$
F_k + (-1)^{n-2} 2 = F_{n+1}L_{n-2},
$$
and hence
\begin{equation} \label{cookie}
F_{n+\delta_1} L_{n-\delta_2}  = y^p \; \mbox{ where } \;
\{ \delta_1, \delta_2 \} = \{ 1, 2 \}.
\end{equation}
We claim that
\begin{equation} \label{fudge}
\gcd ( F_{k+3}, L_{k} ) = 
\left\{
\begin{array}{cc}
4 & \mbox{ if } k \equiv 3 \mod{6}   \\
2 & \mbox{ if } k \equiv 0 \mod{6} \\
1 & \mbox{ otherwise.}
\end{array}
\right.
\end{equation}
To see this, note the identity
$$
3 F_{k+3} = L_k + 4 F_{k+2}
$$
which implies, since $\gcd (F_{k+3},F_{k+2}) =1$,  that $\gcd ( F_{k+3}, L_{k} ) \mid 4$. The fact that
$$
F_{k+3} \equiv L_k \equiv 2 \mod{4} \; \mbox{ if } k \equiv 0 \mod{6}
$$
and
$$
F_{k+3} \equiv L_k \equiv 0 \mod{4} \; \mbox{ if } k \equiv 3 \mod{6},
$$
while  $F_k$ and $L_k$ are odd unless $3 \mid k$
completes the proof.

From (\ref{cookie}) and (\ref{fudge}), it thus follows that $L_{n-\delta_2} = 2^\alpha y_1^p$ for integers $\alpha \ge 0$ and $y_1$.  Appealing to Theorem 2 of \cite{BLMS2}, and the identity $L_{-m}=(-1)^m L_m$, we thus have that 
$$
|n-\delta_2| \in \{ 0, 1, 3, 6 \}.
$$
%and hence
%$$
%n \in \{ -5, -4, -2, -1, 0, 1, 2, 3, 4, 5, 7, 8 \}.
%$$
We check that $F_{2n+1} \pm 2$ is a perfect power only for those $n$ corresponding to
$$
F_{-9}+2 = F_9+2 = 6^2, \; F_{-9}-2=F_9-2 = 2^5, \;  F_{-3}+2 = F_3+2 = 2^2,
$$
$$
F_{-3}-2 = F_3-2 = 0 \; \mbox{ and } \; F_{-1}-2 = F_1-2 = -1.
$$

We may thus suppose that $k=2n$ is even, so that $F_{-k}=-F_{k}$, and hence, without loss of generality, 
that $F_{2n}+2=\pm y^p$. The case $p=2$ is easily dealt with by reducing the problem to the determination
of integral points on elliptic curves. We may therefore suppose $p \ge 3$ and so absorb the sign
into the $y^p$. We therefore consider the equation 
\begin{equation}\label{eqn:main}
F_{2n}+2=y^p.
\end{equation}
This is of the shape (\ref{eqn:ShiftedPowers}) with $k=2n$, $c=2$ and $a=1$. Writing
$x=\varepsilon^{2n}+\sqrt{5}$, equation (\ref{pp2-1}) implies that
\begin{equation}\label{eqn:pp2}
x^2-6=\sqrt{5} \varepsilon^{2n} y^p.
\end{equation}
By thinking of the constant $-6$ as $-6 \cdot 1^p$, we may view this equation
as generalized Fermat equation of signature $(p,p,2)$ over 
$\Q(\sqrt{5})$. To the solution $(x,y,n,p)$ of \eqref{eqn:pp2} (and hence
to the solution $(n,y,p)$ to \eqref{eqn:main}) we associate the Frey--Hellegouarch
equation 
\begin{equation} \label{eqn:Frey}
E_n \; : \; Y^2=X^3+2x X^2+6 X, \qquad x=\varepsilon^{2n}+\sqrt{5}.
\end{equation}
 This will prove much easier to deal with 
than the corresponding $(p,p,p)$ equation defined over $\Q(\sqrt{5},\sqrt{6})$ that we obtain from the arguments of \cite{BDMS}.
We shall apply modularity and level-lowering to mod $p$ representation of $E_n$ to deduce the following.
\begin{prop}\label{prop:ll}
Let $(n,y,p)$ be a solution to \eqref{eqn:main} with $p \ge 5$. Let $\overline{\rho}_{E_n,p}$
be the mod $p$ representation of $E_n$. Then $\overline{\rho}_{E_n,p}$
is irreducible. Moreover, 
%We know that $\overline{\rho}_{E_n,p}$ is
%irreducible and we can apply the usual level-lowering recipe.
$\overline{\rho}_{E_n,p} \sim \overline{\rho}_{\ff,\pi}$
where $\ff$ Hilbert eigenform over $\Q(\sqrt{5})$ of weight $(2,2)$ that is new of level
\begin{equation}\label{eqn:level}
\cN=(2)^7 \cdot (3) \cdot (\sqrt{5});
\end{equation}
here $\pi \mid p$ is some prime of $\OO_\ff$, the ring
of integers of the number field generated by the Hecke eigenvalues
of $\ff$.
\end{prop}
The Hilbert newspace for weight $(2,2)$ and level $\cN$
has dimension $6144$. It is not possible using current
software capabilities to compute the eigenforms belonging 
to this space. One of the novelties of the current paper
is a sieving argument that works with mod $p$ eigensystems
to eliminate all of the space except for three elliptic curves.

\section{Dealing with small $p$ and small $\lvert y \rvert$}
We shall apply the methods of Galois representation
and modularity to the equation \eqref{eqn:main}. Such methods are somewhat harder to apply with small exponent $p$,
and so in this section we deal with the cases $p=2$ and $p=3$ separately. 
Later on we would like to apply bounds for linear forms in logarithms to \eqref{eqn:main}, and for this
it is useful to know that $y$ is not too small. We show below that if 
$n \ne -2$, $-1$ then $\lvert y \rvert \ge 19$.
\begin{lem}\label{lem:peq2}
The only solutions to the equation $F_{2n}+2=\pm y^2$
are $(n,y)=(-1,\pm 1)$ and $(-2,\pm 1)$. 
\end{lem}
\begin{proof}
Let $Y=5 y L_{2n}$ and $X=5 y^2$. It follows from identity \eqref{eqn:wk}
that $(X,Y)$ is an integral point on one of the two elliptic curves
\[
Y^2=X(X^2 - 20 X+120), \qquad Y^2=X(X^2 + 20 X+120).
\]
To determine the integral points on these two elliptic curves
 we used the computer package \texttt{Magma} \cite{magma}
which utilizes a standard algorithm that employs 
lower bounds for linear forms in elliptic logarithms \cite{Smart}.
We
find that the integral points on the first curve are $(0,0)$, $(5, \pm 15)$, $(24, \pm 72)$,
and those on the second are $(0,0)$, $(5,\pm 35)$, $(24,\pm 168)$.
The lemma follows.
\end{proof}

%\begin{lem}
%If $(n,y,p)$ is a solution to \eqref{eqn:main}
%then
%\begin{equation}\label{eqn:Ln}
%L_{2n}^2=5 y^{2p}-20 y^p+24.
%\end{equation}
%\end{lem}
%\begin{proof}
%This follows from the identity \eqref{eqn:wk}.
%\end{proof}

\begin{lem}
If $p=3$ then the only solutions to \eqref{eqn:main}
are $(n,y)=(-1,1)$ and $(n,y)=(-2,-1)$. 
\end{lem}
\begin{proof}
%Equation \eqref{eqn:Ln} leads us to consider integral points
%on the genus $2$ curve $Y^2=5 X^6-20 X^3+24$. The Mordell--Weil
%group of its Jacobian has $2$-Selmer rank $2$ and we are able 
%to find one generator, thus it appears that the rank is in fact $2$,
%making it difficult to find the integral points. We shall in
%fact proceed differently, starting from equation \eqref{eqn:pp2}.
Write $2n=\alpha+3 m$ where $\alpha=0$, $\pm 1$. Let
\[
X=\sqrt{5} \varepsilon^{m+\alpha} y, \qquad Y=\sqrt{5} \varepsilon^\alpha x.
\]
From \eqref{eqn:pp2} we deduce that $(X,Y)$ is an $\OO_K$-integral
point on the elliptic curve
\[
Y^2=X^3+30 \varepsilon^{2 \alpha}.
\]
These three elliptic curves (corresponding to $\alpha=0$, $1$, $-1$)
all have rank $2$ over $K$, and we are able to compute the $\OO_K$-integral
points via an algorithm of Smart and Stephens \cite{SS} implemented
in \texttt{Magma}. These points are 
\[
(19 , \pm 83 ),
\qquad 
((- 3-\sqrt{5})/2 , \pm (1-2\sqrt{5}) ), 
\qquad
((- 3+\sqrt{5})/2 , \pm (1+2\sqrt{5}) ), 
\]
for $\alpha=0$,
and
\begin{align*}
 ((3-5 \sqrt{5})/2 , \pm( 8-5 \sqrt{5}) ), \qquad & 
(\sqrt{5} , \pm (5+2 \sqrt{5}) ), \\ 
((5+3 \sqrt{5})/2 , \pm ( 10+3 \sqrt{5} )), \qquad & 
((55+15 \sqrt{5})/2 , \pm (165+58 \sqrt{5}) ) \, ,
\end{align*}
for $\alpha=1$, with conjugate points for $\alpha=-1$.
The lemma easily follows.
\end{proof}
\begin{lem}\label{lem:div}
Suppose $(n,y,p)$ is a solution to \eqref{eqn:main}.
If $q \mid y$
is prime, then $q \equiv 1$, $5$, $19$, $23 \pmod{24}$.
In particular,
$2 \nmid y$ and $3 \nmid y$. Moreover, 
\[
n \equiv 2,\, 4,\, 7,\, 8,\, 10,\, 11 \pmod{12}
\]
\end{lem}
\begin{proof}
Suppose $2 \mid y$. From $F_{2n}+2=y^p$ we have $2 \mid\mid F_{2n}$.
However, $2 \mid F_{2n}$ implies that $3 \mid n$. Thus 
$F_6 \mid F_{2n}$. As $F_6=8$
we have a contradiction.

Now suppose $q \mid y$ is an odd prime. From \eqref{eqn:wk} we obtain
\[
L_{2n}^2=5 y^{2p}-20 y^p+24,
\]
so $24 \equiv L_{2n}^2 \pmod{q^2}$. Thus $q \ne 3$,
and $6$ is a quadratic residue modulo $q$. It follows
that $q \equiv 1$, $5$, $19$, $23 \pmod{24}$.

The final part of the lemma follows from considering
$F_{2n}+2$ modulo $6$. 
\end{proof}

\begin{lem}\label{lem:powerof5}
The only solutions to the equation $F_{2n}+2=\pm 5^m$ are
$F_{-4}+2=-1$, $F_{-2}+2=1$, $F_{4}+2=5$.
%$(n,m)=(-1,0)$, $(2,1)$.
\end{lem}
\begin{proof}
As above we deduce that
\[
L_{2n}^2=5 \cdot 5^{2m} \mp 20 \cdot 5^m+24.
\]
If $m$ is even then write
\[
X=5^{m+1}, \qquad Y=5^{(m+2)/2} \cdot L_{2n}.
\]
Then $(X,Y)$ satisfies 
\[
Y^2=X^3 \mp 20 X^2+120X;
\]
we are interested in computing the integral points on these two elliptic
curves. For this we used the computer package \texttt{Magma} \cite{magma}
which utilizes a standard algorithm to determine integral points via
lower bounds for linear forms in elliptic logarithms \cite{Smart}.
The integral points on the model 
$Y^2=X^3-20X^2+120X$ are $(0,0)$, $(5,\pm 15)$, $(24,\pm 72)$,
and lead to the solution $F_{-2}+2=1$. 
The integral points on the model 
$Y^2=X^3+20X^2+120X$ are $(0,0)$, $(5,\pm 35)$, $(24,\pm 168)$,
and lead to the solution $F_{-4}+2=-1$. 

If $m$ is odd then write
\[
X=5^{m+2}, \qquad Y=5^{(m+5)/2} \cdot L_{2n}.
\]
Then $(X,Y)$ satisfies 
\[
Y^2=X^3 \mp 100 X^2+3000X.
\]
The integral points on the model 
$Y^2=X^3 - 100 X^2+3000X$
are $(0,0)$, $(24,\pm 168)$, $(125,\pm 875)$
and lead to the solution $F_{4}+2=5$.
The integral points on the model 
$Y^2=X^3 + 100 X^2+3000X$
are
$(0 , 0)$, $(2904 , \pm 159192)$ and do not lead to any
solutions to the original equation.
\end{proof}

%\begin{lem}\label{lem:smallsolutions}
%The only solutions to \eqref{eqn:main}
%with $\lvert n \vert \le 1000$ are
%$F_{-2}+2=1$ and $F_{-4}+2=-1$.
%\end{lem}
%\begin{proof}
%We checked this by a direct computation using \texttt{Magma} \cite{magma}.
%\end{proof}

\begin{lem}\label{lem:yge19}
Let $(n,y,p)$ be a solution to \eqref{eqn:main} and
suppose that $n \ne -2$, $-1$. Then $\lvert y \rvert \ge 19$.
\end{lem}
\begin{proof}
As $n \ne -2$, $-1$ it follows that $\lvert y \rvert>1$.
Suppose $\lvert y \rvert<19$. By Lemma~\ref{lem:div}, the
only prime divisor of $y$ is $5$. This now contradicts Lemma~\ref{lem:powerof5}.
\end{proof}

%---------------------------------------
%\section{Preparation for Sieving}\label{sec:sieveprep}
%---------------------------------------

%--------------------------------------------
\section{Irreducibility of the mod $p$ representation}
%-------------------------------------------

Henceforth $(n,y,p)$ is a solution to \eqref{eqn:main}
with prime exponent $p \ge 5$, and $E_n$
is the Frey--Hellegouarch curve $E_n$ given by \eqref{eqn:Frey}.
An easy application of Tate's algorithm (together with Lemma~\ref{lem:div}) 
yields the following.
\begin{lem}\label{lem:Tate}
The model in \eqref{eqn:Frey} is minimal with discriminant
and conductor
\[
\Delta=2^8 \cdot 3^2 \cdot \varepsilon^{2n} \cdot \sqrt{5} \cdot y^p,
\qquad
\mathfrak{N}=(2)^7 \cdot (3) \cdot (\sqrt{5}) \cdot 
\prod_{\fq \mid y,\, \fq \ne (\sqrt{5})}
\fq.
\]
\end{lem}
We would like to apply level-lowering to the mod $p$ representation $\overline{\rho}_{E_n,p}$,
and for this we need to show that it is irreducible. 
We shall make use of the following result
due to Freitas and Siksek \cite{FS2}, 
which is based on the work of David \cite{DavidI}
and Momose \cite{Momose}.
\begin{prop}\label{prop:irred}
Let $K$ be a totally real Galois number field of degree $d$,
with ring of integers $\OO_K$ and Galois group $G=\Gal(K/\Q)$.
Let $S=\{0,12\}^G$, which we think of as the set of sequences of values $0$, $12$
indexed by $\tau \in G$. %Fix an ordering $\tau_1,\tau_2,\ldots,\tau_d$
For $\bs=(s_\tau) \in S$ and $\alpha \in K$, define the \textbf{twisted norm associated
to $\bs$} by
\[
\cN_\bs(\alpha)= \prod_{\tau \in G} \tau(\alpha)^{s_\tau}.
\]
Let $\varepsilon_1,\dots,\varepsilon_{d-1}$
be a basis for the unit group of $K$,
and define
\begin{equation}\label{eqn:As}
A_\bs:=\Norm  \left( \gcd ( ( \cN_\bs(\varepsilon_1)-1) \OO_K,\ldots, (\cN_\bs(\varepsilon_{d-1})-1  ) \OO_K) \right).
\end{equation}
Let $B$ be the least common multiple of the $A_\bs$ taken over all $\bs \ne (0)_{\tau \in G}$,
$(12)_{\tau \in G}$.
Let $p \nmid B$ be a rational prime, unramified in $K$, such that $p \geq 17$ or $p = 11$.
Let $E/K$ be an elliptic curve, and $\fq \nmid p$ be a prime of good reduction for $E$.
Define
\[
P_\fq(X)=X^2-a_\fq(E) X + \Norm(\fq)
\]
to be the characteristic polynomial
of Frobenius for $E$ at $\fq$. Let $r \ge 1$ be an integer such that
$\fq^r$ is principal.
If $E$ is semistable at all $\fp \mid p$
and $\overline{\rho}_{E,p}$ is reducible then
\begin{equation}\label{eqn:res}
p \mid \Res(\, P_\fq(X) \, , \, X^{12 r}-1\, )
\end{equation}
where $\Res$ denotes the resultant of the two polynomials.
\end{prop}
We observe in passing that since $P_\fq(X)$ has two complex roots
of absolute value $\sqrt{\Norm(\fq)}$, the resultant in \eqref{eqn:res}
cannot be zero. We now arrive at the main result of this section.
\begin{lem}\label{lem:irred}
Let $(n,y,p)$ be a solution to \eqref{eqn:main}
with $p \ge 5$ prime.
Let $E_n$
be the Frey curve given in \eqref{eqn:Frey}.
Then $\overline{\rho}_{E_n,p}$ is irreducible.
\end{lem}
\begin{proof}
Let 
\begin{equation}\label{eqn:M0}
M_0=2520=2^3 \cdot 3^2 \cdot 5 \cdot 7 \, ,
\end{equation}
and
\[
\cQ^\prime=\{ \text{$\fq$ a prime ideal of $\OO_K$} \quad : \quad \text{$\fq \nmid 2 \cdot 3 \cdot \sqrt{5}$} \quad \text{ and} \quad \text{ $\norm(\fq)<300$}  \}.
\]
Let
\begin{equation}\label{eqn:cQ}
\cQ=\{\fq \in \cQ^\prime \quad : \quad
\text{the multiplicative order 
of $\varepsilon^2$
in $\F_\fq$ divides $M_0$}\}.
\end{equation}
The set $\cQ$ contains $25$ prime ideals $\fq$. The Frey curve 
\eqref{eqn:Frey}
modulo $\fq$ depends only on $n$ modulo $M_0$.
Let
\begin{equation}\label{eqn:cM}
\cM=\{ 0 \le m \le M_0 \; : \; m \equiv 2,4,7,8,10,11 \pmod{12}  \}.
\end{equation}
By Lemma~\ref{lem:div}, if $(n,y,p)$ is a solution to
\eqref{eqn:main} then $n \equiv m \pmod{M_0}$ for some
unique $m \in \cM$. 
%We prove the lemma first for  $p=5$, $7$, $13$. Then $n \equiv m
%\pmod{M_0}$ for some (unique) $m \in \cM$. 
%Let $\fq \in \cQ$.
In particular, $\varepsilon^{2n} \equiv \varepsilon^{2m} \pmod{\fq}$.
Suppose 
$\fq \nmid ((\varepsilon^{2m}+\sqrt{5})^2-6)$. 
By \eqref{eqn:pp2}, we have $\fq \nmid y$. By Lemma~\ref{lem:Tate}
we see that
$E_n$ has good reduction modulo $\fq$. Moreover, $a_\fq(E_n)=a_\fq(E_m)$.
In particular, if $t^2-a_\fq(E_m)t+\norm(\fq)$ is irreducible
modulo $p$, then $\overline{\rho}_{E_n,p}$ is irreducible.

We wrote a short \texttt{Magma} script which did the following.
For each of the values $p=5$, $7$, $13$, and for each $m \in \cM$,
it verified that there an $\fq \in \cQ$ such that
$\fq \nmid ((\varepsilon^{2m}+\sqrt{5})^2-6)$
and that $t^2-a_\fq(E_m)t+\norm(\fq)$ is irreducible
modulo $p$. This completes the proof for $p=5$, $7$, $13$.

Thus we suppose that $p=11$ or $p \ge 17$.
We apply the above proposition. A fundamental unit for
$K=\Q(\sqrt{5})$ is $\varepsilon$, and it follows that $B=320$,
where $B$ is as in the statement of the proposition.
Thus $p \nmid B$. Moreover, from Lemma~\ref{lem:Tate},
$E_n$ is semistable at $\fp \mid p$. We suppose that
$\overline{\rho}_{E_n,p}$ is reducible.
Let 
\[
%S=\{\fq_1,\; \fq_2,\; \fq_5,\; \fq_6,\; \fq_7,\; \fq_8,\;
% \fq_{10},\; \fq_{11}\},
S=\{\fq \in \cQ \; : \; 
\text{$\fq$ is above a rational prime $q \not\equiv 1,5,19,23$}\}.
\]
The set $S$ has size $15$.
%where $\fq_i$ are as in the previous section. Here
%we have chosen the $\fq_i$ that are above a rational
%prime $q \not\equiv 1$, $5$, $19$, $23 \pmod{24}$.
By Lemma~\ref{lem:div}, it follows that %$\fq \nmid y$
%and 
$E_n$ has good reduction at all $\fq \in S$.
Recall that $n \equiv m \pmod{M_0}$ for some unique $m \in \cM$.
Moreover, $a_\fq(E_n)=a_\fq(E_m)$ for $\fq \in S$.
It follows from the above proposition that $p$ divides
\[
\gcd(\{ 
\Res(t^2-a_\fq(E_m)t+\norm(\fq),t^{12}-1) \; : \;
\fq \in S\}).
\]
We computed this greatest common divisor for each $m \in \cM$
and verified that it is never divisible by $11$ or any prime $\ge 17$.
The lemma follows.
\end{proof}

\section{Level-lowering and consequences}

We are now in a position to prove Proposition~\ref{prop:ll}
%We assume (for now) that the mod $p$ representation for $E_n$ is
\begin{proof}[Proof of Proposition~\ref{prop:ll}]
The elliptic curve $E_n$ is modular by \cite{FLS},
and mod $p$ representation $\overline{\rho}_{E_n,p}$ 
is irreducible by Lemma~\ref{lem:irred}.
If $p>5$ then the proposition immediately follows from the statement of
Theorem 7 of \cite{FS} 
which is based on the works of Fujiwara, Jarvis and Rajaei.
Now let $p=5$. In this case the statement of theorem in \cite{FS}
is inapplicable in our situation. Specifically condition (v)
of that theorem is not satisfied in our setting as
$5 \nmid \ord_{\sqrt{5}}(\Delta)$. However that condition 
is only needed to remove the primes above $p$ from the level without
increasing the weight. In our situation we content ourselves, when $p=5$,
with removing from the level the primes dividing $y$ that do not also divide
$\nmid 2\cdot 3 \cdot \sqrt{5}$. As in \cite{FS} this can be done
whilst keeping the weight $(2,2)$.
\end{proof}
\begin{lem}\label{lem:congs}
With notation as in Proposition~\ref{prop:ll}, let $\fq \nmid p \cdot \cN$ be a prime of $\OO_K$.
Let $m$ be an integer satisfying $\varepsilon^{2m} \equiv \varepsilon^{2n} \pmod{\fq}$.
Write
\begin{equation}\label{eqn:bm}
b_\fq(m)=
\begin{cases}
a_\fq(E_m) & \text{$\fq \nmid (F_{2m}+2)$}\\
\norm(\fq)+1 & \text{$\fq \mid (F_{2m}+2)$  and $-(\varepsilon^{2m}+\sqrt{5})$}\\
& \text{is a square modulo $\fq$} \\
-\norm(\fq)-1 & \text{otherwise}.
\end{cases}
\end{equation}
\begin{comment}
\begin{itemize}
\item 
%If $\fq \nmid 3 \cdot \sqrt{5} \cdot (F_{2n}+2)$ then 
If $\fq \nmid (F_{2n}+2)$ then 
\[
a_\fq(E_n) \equiv a_\fq(\ff) \pmod{\pi}.
\]
%\item If $\fq \mid 3 \cdot \sqrt{5} \cdot (F_{2n}+2)$ and $-(\varepsilon^{2n}+\sqrt{5})$
\item If $\fq \mid (F_{2n}+2)$ and $-(\varepsilon^{2n}+\sqrt{5})$
is a square modulo $\fq$ then
\[
\Norm(\fq)+1 \equiv a_\fq(\ff) \pmod{\pi}.
\]
%\item If $\fq \mid 3 \cdot \sqrt{5} \cdot (F_{2n}+2)$ and $-(\varepsilon^{2n}+\sqrt{5})$
\item If $\fq \mid  (F_{2n}+2)$ and $-(\varepsilon^{2n}+\sqrt{5})$
is a non-square modulo $\fq$ then
\[
-\Norm(\fq)-1 \equiv a_\fq(\ff) \pmod{\pi}.
\]
\end{itemize}
\end{comment}
Then $b_\fq(m) \equiv a_\fq(\ff) \pmod{\pi}$.
\end{lem}
\begin{proof}
Suppose $\fq \nmid p\cdot \cN$.
Since $F_{2n}+2=y^p$, we see from Lemma~\ref{lem:Tate}
that $E_n$ has good reduction at $\fq$ if 
%$\fq \nmid 3 \cdot \sqrt{5} \cdot (F_{2n}+2)$
$\fq \nmid (F_{2n}+2)$
and multiplicative reduction at $\fq$ if 
%$\fq \mid 3 \cdot \sqrt{5} \cdot 
$\fq \mid (F_{2n}+2)$. 
%In the former case we obtain (i) by taking traces of the
%images of Frobenius at $\fq$ in 
%$\overline{\rho}_{E_n,p} \sim \overline{\rho}_{E,p}$.
Suppose we are in the latter case.
We know \cite[Theorem V.5.3]{SilvermanII} that the
reduction at $\fq$ is split if and only if $-c_6/c_4$
is a $\fq$-adic square, where $c_4$ and $c_6$
are the usual $c$-invariants of $E_n$. In our case
\[
c_4=2^5(2 x^2-9), \qquad c_6=2^7 x (-4 x^2+27).
\]
%Suppose first that $\fq \ne 3 \OO_K$. Thus 
From \eqref{eqn:pp2}
we have $x^2 \equiv 6 \pmod{\fq}$ and so $-c_6/c_4 \equiv -x \pmod{\fq}$.
As $x=\varepsilon^{2n}+\sqrt{5}$, the multiplicative reduction
at $\fq$ is split if and only if $-(\varepsilon^{2n}+\sqrt{5})$
is a square modulo $\fq$. 
%Suppose now that $\fq=3 \OO_K$.
%From \eqref{eqn:pp2} we see that  $\fq \nmid$ and so
%$-c_6/c_4 \equiv 8x \equiv -x \pmod{\fq}$ giving
%the same conclusion as before.

By comparing the traces of the images
of the Frobenius element at $\fq$ in
$\overline{\rho}_{E_n,p} \sim \overline{\rho}_{\ff,\pi}$
we obtain $b_\fq(n) \equiv a_\fq(\ff) \pmod{\pi}$
in all cases. Finally, as $\varepsilon^{2m} \equiv \varepsilon^{2n} \pmod{\fq}$,
it follows that $F_{2m} \equiv F_{2n} \pmod{\fq}$, and so $b_\fq(m)=b_\fq(n)$
proving the lemma.
\end{proof}

Write $S=S_{(2,2)}^{\mathrm{new}}(\cN)$. Using
\texttt{Magma}
we find that $S$ has dimension $6144$. 
We let $\cF$ be the set of eigenforms $\ff$ belonging to $S$ 
(thus $\#\cF=6144$).
Alas it is not practical to compute these newforms with
current software capabilities. However it is quite practical
using \texttt{Magma} 
to compute the action of the Hecke operators $T_\fq$ on $S$
for small primes $\fq$ of $\OO_K$. For the theoretical details
behind these algorithms we recommend \cite{DV}.

We used a \texttt{Magma} program written by Stephen Donnelly
to search for elliptic curves over number fields
with a given conductor.
This program found $288$ pairwise 
non-isogenous elliptic curves $F/K$ with conductor $\cN$.
We know by \cite{FLS} that these corresponds to
$288$ distinct $\ff \in \cF$ with rational Hecke eigenvalues.
We let $\cE$ be this set of these $288$ elliptic curves and we let 
$\cF^\prime$ be the subset of $\cF$ coming from these
$288$ elliptic curves.

\section{Reducing to elliptic curves}
\begin{prop}\label{prop:rat}
Let $(n,y,p)$ be a solution to \eqref{eqn:main} with prime exponent
$p \ge 5$. Then $\overline{\rho}_{E_n,p} \sim \overline{\rho}_{E,p}$
where $E \in \cE$. 
\end{prop}
We shall prove Proposition~\ref{prop:rat} by contradiction. 
Suppose $\overline{\rho}_{E_n,p} \not \sim \overline{\rho}_{E,p}$
for any $E \in \cE$. Then 
$\overline{\rho}_{E_n,p} \sim \overline{\rho}_{\fg,\pi}$
for some $\fg \in \cF - \cF^\prime$. 
Let $\cQ$ and $\cM$ be as in the proof of Lemma~\ref{lem:irred}.
Let $m$ be the unique
element of $\cM$ such that $m \equiv n \pmod{M_0}$.
In particular, we know that
$\varepsilon^{2m} \equiv \varepsilon^{2n} \pmod{\fq}$ for all $\fq \in \cQ$.
From Lemma~\ref{lem:congs} we see that
\begin{equation}\label{eqn:bmcong}
b_\fq(m) \equiv a_\fq(\fg) \pmod{\pi}
\end{equation}
for all $\fq \in \cQ$ with $\fq \nmid p$. 

Suppose for now that $\fq \in \cQ$ and $\fq \nmid p$. Write $T_\fq$ for the Hecke operator
corresponding to $\fq$ acting on the space $S=S_{(2,2)}^{\mathrm{new}}(\cN)$.
Let $C_\fq(x)=\det(x I - T_\fq) \in \Z[x]$ be its characteristic polynomial; this has
roots $a_\fq(\ff)$ with $\ff$ running through $\ff \in \cF$. Now let
\[
C^\prime_\fq(x)=\prod_{E \in \cE} (x-a_\fq(E)) \in \Z[x].
\]
Thus $C^\prime_\fq(x)$ divides $C_\fq(x)$. Moreover, let
\[
C^{\prime\prime}_\fq(x)=\frac{C_\fq(x)}{C^\prime_\fq(x)} \in \Z[x].
\]
The roots of $C^{\prime\prime}_\fq(x)$ are $a_\fq(\ff)$ with $\ff$ running through $\ff \in \cF^\prime$.
We see from \eqref{eqn:bmcong} that $C^{\prime\prime}(b_\fq(m)) \equiv 0 \pmod{\pi}$.
However as $C^{\prime\prime} \in \Z[x]$ and $b_\fq(m) \in \Z$ it follows that 
\[
C^{\prime\prime}(b_\fq(m)) \equiv 0 \pmod{p}.
\]
Now let 
\[
G_{m,\fq}=\norm(\fq) \cdot C^{\prime\prime}(b_\fq(m)) \in \Z.
\]
We see that $p \mid G_{m,\fq}$ for all $\fq \in \cQ$ regardless of whether $\fq$ divides $p$ or not.
Thus $p$ divides
\[
H_m:=\gcd \{ G_{m,\fq} \; : \; \fq \in \cQ \}.
\]
We computed the integers $H_m$ for all $m \in \cM$ and factorized them. It turns that all are non-zero,
which means we have bounded $p$ under the assumption that $\overline{\rho}_{E_n,p} \not \sim \overline{\rho}_{E,p}$
for all $E \in \cE$. In particular it turns out $p \le 109$. More precisely, we are left to consider
$9391$ pairs $(p,m)$ where $p\ge 5$ is a prime dividing $H_m$.  

To proceed further we remark that Hilbert Modular Forms package in 
\texttt{Magma} computes a matrix, which shall denote by $R_\fq$, giving the action of
the operator $T_\fq$ (with $\fq$ not dividing the level $\cN$)
with respect to a $\Z$-basis of a lattice in $S_{(2,2)}^{\mathrm{new}}(\cN)$ that is Hecke-stable.
Write $\overline{R}_\fq$ for the reduction of $R_\fq$ modulo $p$. Write $\overline{\fg}$ for the
mod $p$ eigensystem corresponding to $\fg$. It follows from the above that the intersection
\begin{equation}\label{eqn:intersection}
\bigcap_{\substack{{\fq \in \cQ}\\{\fq \nmid p}}} \Ker \left(\overline{R}_{\fq}-\overline{b_\fq(m)}\cdot I\right)
\end{equation}
contains an $\F_p$-line corresponding to $\overline{\fg}$. 
We computed the intersection \eqref{eqn:intersection}
for all the $9391$ remaining pairs $(p,m)$. This is merely $\F_p$-linear algebra once the matrices $R_\fq$
representing the Hecke operators were computed. We found that for all but $21$ of the $9391$ pairs
$(p,m)$ that the space \eqref{eqn:intersection} is $0$-dimensional giving us a contradiction.
For the proof of Proposition~\ref{prop:rat} we need now only consider the following $21$ remaining pairs $(p,m)$:
\begin{gather*}
(5,2),\quad (5, 2518),\quad (5, 2519),\quad (7,2),\quad (7, 2518), \quad (7, 2519),\\
(11,2), \quad (11, 2518), \quad (11, 2519),\quad (13,2),\quad (13, 2518), 
\quad (13, 2519),\\
(17,2),\quad (17, 2519), \quad (19,2),\quad (19, 2519), \quad  (23, 2518), \quad (29, 2518),\\
 (29, 2519), \quad
(41,2),\quad (43, 2518). 
\end{gather*} 
Observe that 
$m \equiv 2$, $-2$, $-1 \pmod{M_0}$ in every one of these $21$ cases.
%The values of $m$ in these $21$ cases are $2$, $2518 \equiv -2 \pmod{M_0}$ and $2519 \equiv -1 \pmod{M_0}$.
The presence of the possibilities $-2$, $-1$ is hardly surprising in view of the solutions $(n,y,p)=(-2,-1,p)$ and $(-1,1,p)$
to \eqref{eqn:main}; for an explanation of the value $2$ see the next
section. 
In all these $21$ cases we found that the intersection \eqref{eqn:intersection}
is $1$-dimensional. We let $E$ be $E_2$ if $m=2$,
$E_{-2}$ if $m=2518 \equiv -2 \pmod{M_0}$ and
$E_{-1}$ if $m=2519 \equiv -1 \pmod{M_0}$. These all have conductors 
$\cN$. Let $\ff \in \cF$ be the Hilbert eigenform
corresponding to $E$. Then $b_\fq(m) \equiv a_\fq(E)
=a_\fq(\ff) \pmod{p}$. It follows that the reduction of the line
corresponding to $\ff$ belongs to the $1$-dimensional 
intersection \eqref{eqn:intersection}, which also contains
the reduction of the line corresponding to $\fg$. 
Thus the mod $p$ eigensystems $\overline{\ff}$ and $\overline{\fg}$
are equal. It follows that 
%the mod $p$ Galois representations
$\overline{\rho}_{\ff,p}\sim \overline{\rho}_{\fg,\pi}$.
Thus $\overline{\rho}_{E_n,p} \sim \overline{\rho}_{E,p}$.
But $E_2$, $E_{-2}$, $E_{-1} \in \cE$; this completes the proof of Proposition~\ref{prop:rat}.

\section{Reducing to only three elliptic curves}
We know from Proposition~\ref{prop:rat} that 
$\overline{\rho}_{E_n,p} \sim \overline{\rho}_{E,p}$
where $E$ is one of the $288$ elliptic curves belonging to $\cE$.
In this section we eliminate all but three of the elliptic curves
belonging to $\cE$. 
\begin{prop}\label{prop:threeells}
Let $(n,y,p)$ be a solution to \eqref{eqn:main} with $p \ge 5$.
Then $n \equiv m \pmod{M_0}$ and 
$\overline{\rho}_{E_n,p} \sim \overline{\rho}_{E,p}$ where
\begin{enumerate}
\item[(i)] $m=2$ and $E=E_2$;
\item[(ii)] $m=M_0-2$ and $E=E_{-2}$;
\item[(iii)] $m=M_0-1$ and $E=E_{-1}$.
\end{enumerate}
\end{prop}
We shall need the following slight
strengthening of Lemma~\ref{lem:congs}.
\begin{lem}\label{lem:congs2}
Let $(n,y,p)$ be a solution to \eqref{eqn:main} with $p \ge 5$.
Let $E \in \cE$ satisfy $\overline{\rho}_{E_n,p} \sim \overline{\rho}_{E,p}$.
Let $\fq \nmid \cN$ be a prime of $\OO_K$. Let $m$ be an integer
satisfying $\varepsilon^{2m} \equiv \varepsilon^{2n} \pmod{\fq}$.
Then $b_\fq(m) \equiv a_\fq(E) \pmod{p}$,
where $b_\fq(n)$ is given by \eqref{eqn:bm}.
\begin{comment}
\begin{itemize}
\item 
%If $\fq \nmid 3 \cdot \sqrt{5} \cdot (F_{2n}+2)$ then 
If $\fq \nmid (F_{2n}+2)$ then 
\[
a_\fq(E_n) \equiv a_\fq(E) \pmod{p}.
\]
%\item If $\fq \mid 3 \cdot \sqrt{5} \cdot (F_{2n}+2)$ and $-(\varepsilon^{2n}+\sqrt{5})$
\item If $\fq \mid (F_{2n}+2)$ and $-(\varepsilon^{2n}+\sqrt{5})$
is a square modulo $\fq$ then
\[
\Norm(\fq)+1 \equiv a_\fq(E) \pmod{p}.
\]
%\item If $\fq \mid 3 \cdot \sqrt{5} \cdot (F_{2n}+2)$ and $-(\varepsilon^{2n}+\sqrt{5})$
\item If $\fq \mid  (F_{2n}+2)$ and $-(\varepsilon^{2n}+\sqrt{5})$
is a non-square modulo $\fq$ then
\[
-\Norm(\fq)-1 \equiv a_\fq(E) \pmod{p}.
\]
\end{itemize}
\end{comment}
\end{lem}
\begin{proof}
If $\fq \nmid p$ then this is a special case of Lemma~\ref{lem:congs}. If $\fq \mid p$
then this follows from the proof of Lemma~\ref{lem:congs} together with \cite{KO}.
\end{proof}
Now let $\cQ$ be as in the proof of Lemma~\ref{lem:irred}. The following is immediate.
\begin{lem}\label{lem:BmE}
Let $(n,y,p)$ be a solution to \eqref{eqn:main} with $p \ge 5$
prime. Let $E \in \cE$ such that $\overline{\rho}_{E_n,p} \sim \overline{\rho}_{E,p}$.
Let $n \equiv m \pmod{M_0}$ with $m \in \cM$ then
$p$ divides
\[
B_m(E)=\gcd(\{ b_\fq(m)-a_\fq(E) \; : \; \fq \in \cQ \}).
\]
\end{lem}

We computed $B_m(E)$ for all of the $288$ elliptic curves
$E \in \cE$ and $m \in \cM$. We found that $B_m(E)$ is not divisible
by any primes $p \ge 5$ except in three cases
where $B_m(E)=0$:
\begin{enumerate}
\item[(i)] $m=2$ and $E=E_2$;
\item[(ii)] $m=M_0-2$ and $E=E_{-2}$;
\item[(iii)] $m=M_0-1$ and $E=E_{-1}$.
\end{enumerate}
The possibilities (ii) and (iii) are natural, and they
correspond to the solutions $(n,y,p)=(-2,-1,p)$
and $(-1,1,p)$ respectively. The possibility (i)
results from $F_4+2=5$ from which it is easy to deduce
that $E_2$ has conductor $\cN$ and so it is natural
(though annoying) that our sieve cannot eliminate this possibility.
This proves Proposition~\ref{prop:threeells}.

\section{Enlarging $M_0$}
We let 
\begin{equation}\label{eqn:M1}
M_1=M_0 \; \times
\prod_{\substack{\text{$\ell$ prime}\\{11 \le \ell < 10^4}}} \ell ,
\end{equation}
where $M_0$ is given in \eqref{eqn:M0}.
In this section we prove the following.
\begin{lem}\label{lem:M1}
Let $(n,y,p)$ be a solution to \eqref{eqn:main} with $p \ge 5$
prime. Then $\overline{\rho}_{E_n,p} \sim \overline{\rho}_{E,p}$
and $n \equiv m_0 \pmod{M_1}$ where
\begin{enumerate}
\item[(i)] either $m_0=2$ and $E=E_2$;
\item[(ii)] or $m_0=-2$ and $E=E_{-2}$;
\item[(iii)] or $m_0=-1$ and $E=E_{-1}$.
\end{enumerate}
\end{lem}
\begin{proof}
Fix $m_0 \in \{-2,-1,2\}$, let $E=E_{m_0}$
and suppose $\overline{\rho}_{E_n,p} \sim \overline{\rho}_{E,p}$.
We would like to show that $n \equiv m_0 \pmod{M_1}$.

There are $164$ primes in the interval $[11,10000]$; we denote them
by
\[
\ell_1=11, \; \ell_2=13,\; \dotsc,\; \ell_{164}=9973.
\]
We let $L_0=M_0$, and $L_{i}=\ell_i \cdot L_{i-1}$ for $1 \le i \le 164$.
Then $L_{164}=M_1$. 
We shall show inductively
that $n \equiv m_0 \pmod{L_i}$ for $0 \le i \le 164$
which gives the lemma. 
We know by the previous section that
$n \equiv m_0 \pmod{L_0}$. 
For the inductive
step, suppose $n \equiv m_0 \pmod{L_{i-1}}$ and we want to
show that $n \equiv m_0 \pmod{L_{i}}$. Let $\cQ_i$ 
be a set of prime ideals $\fq \nmid \cN$ satisfying the following
\begin{enumerate}
\item[(i)] $\norm(\fq)=q$ is a rational prime $\equiv 1 \pmod{5}$;
\item[(ii)] $\ell_i \mid (q-1)$ and $(q-1) \mid L_{i}$.
%\item[(iii)] the multiplicative order of $\varepsilon^2$ modulo $\fq$
%is divisible by $\ell$.
\end{enumerate}
Let 
\[
\cM_i=\{ 0 \le m \le L_i-1 \; : \;  m \equiv m_0 \pmod{L_{i-1}}\}.
\]
Thus $n \equiv m \pmod{L_i}$ for some unique $m \in  \cM_i$.
Moreover, it follows from (i) and (ii) that $\varepsilon^{2n} \equiv \varepsilon^{2m} \pmod{\fq}$
for all $\fq \in \cQ_i$. 
Define
\[
B_m(\cQ_i) := \gcd\{b_\fq(m)-a_\fq(E) \; : \;  \fq \in \cQ_i\}.
\]
By Lemma~\ref{lem:congs2}, $p \mid B_m(\cQ_i)$.
%Adapting the argument that led us to Lemma~\ref{lem:BmE}
%we see that if $n \equiv m \pmod{L_i}$ then $p \mid B_m$.
We wrote a simple \texttt{Magma} script which for each $1 \le i \le 164$
and for each $m_0 \in \{-2,-1,2\}$ found a set $\cQ_i$ satisfying (i), (ii),
such that, for all $m \in \cM_i$ with $m \not\equiv m_0 \pmod{L_i}$,
the integer $B_m(\cQ_i)$ is non-zero and divisible only by the primes $2$, $3$.
Our computation took a total of around 45 minutes.
This proves the
inductive step and completes the proof.
\end{proof}

\begin{comment}
Before we prove the lemma we need to set up some notation.
If $S_1 \subset \Z/N_1 \Z$ and $S_2 \subset \Z/N_2 \Z$
we define
\[
S_1 \cap S_2  = \{ n \in \Z/N\Z \; : \; (\text{$n$ mod $N_1$}) \in S_1,
\; (\text{$n$ mod $N_2$}) \in S_2  \},
\]
where $N=\lcm(N_1,N_2)$.
\end{comment}

\section{Linear Forms in Three Logs}\label{sec:threelogs}

\begin{comment}
For this section it might help to have a lower bound for $\lvert y \rvert$ in terms of $p$.
\begin{lem}
Let $(n,y,p)$ be a solution to \eqref{eqn:main} with $p \ge 5$. Then
either $y=\pm 1$ or $\lvert y \rvert \ge \sqrt{p}-1$.
\end{lem}
\begin{proof}
Suppose $y \ne \pm 1$. From Lemma~\ref{lem:powerof5} we know that $\lvert y \rvert$ is not a power of $5$.
Let $q$ be a rational prime dividing $y$. From Lemma~\ref{lem:div} we see that $q \ge 19$. Let $\fq$
be a prime of $\OO_K$ above $q$. Thus $\fq \nmid \cN$.
Applying Lemma~\ref{lem:congs2} we have $\norm(\fq)+1 \equiv \pm a_\fq(E) \pmod{p}$
where $E$ is one of three elliptic curves $E_{-2}$, $E_{-1}$, $E_2$.
As $\lvert a_\fq(E) \vert \le 2\sqrt{\norm(\fq)}$ we have
\[
p \le \norm(\fq)+2\sqrt{\norm(\fq)}+1=(\sqrt{\norm(\fq)}+1)^2.
\]
The lemma follows as $\norm(\fq)=q$ or $q^2$ and $q \mid y$.
\end{proof}
\end{comment}
For any algebraic number $\alpha$ of degree $d$ over $\mathbb{Q}$, we define the {\it absolute logarithmic height} of $\alpha$ via the formula
\begin{equation}\label{eqn:htdef}
h(\alpha)= \dfrac{1}{d} \left( \log \vert a_{0} \vert + \sum\limits_{i=1}^{d} \log \max \left( 1, \vert \alpha^{(i)}\vert \right) \right), 
\end{equation}
where $a_{0}$ is the leading coefficient of the minimal polynomial of $\alpha$ over $\mathbb{Z}$ and the $\alpha^{(i)}$ are the conjugates of $\alpha$ in $\mathbb{C}$.
The following is the main result (Theorem 2.1) of Matveev \cite{Mat}.
\begin{thm}[Matveev] \label{Matveev} 
Let $\mathbb{K}$ be an algebraic number field of degree $D$ over $\mathbb{Q}$ and put $\chi=1$ if $\mathbb{K}$ is real, $\chi=2$ otherwise. Suppose that $\alpha_1, \alpha_2, \ldots, \alpha_{n_0} \in \mathbb{K}^*$ with absolute logarithmic heights $h(\alpha_i)$ for $1 \leq i \leq n_0$, and suppose that
$$
A_i \geq \max \{ D \, h (\alpha_i), \left| \log \alpha_i \right| \}, \; 1 \leq i \leq n_0,
$$
for some fixed choice of the logarithm. Define
$$
\Lambda = b_1 \log \alpha_1 + \cdots + b_{n_0} \log \alpha_{n_0},
$$
where the $b_i$ are integers and  set
$$
B = \max \{ 1, \max \{ |b_i| A_i/A_{n_0} \; : \; 1 \leq i \leq n_0 \} \}.
$$
Define, with $e := \exp(1)$, further, 
$$
\Omega =A_1 \cdots A_{n_0}, 
$$
$$
C(n_0) = C(n_0,\chi) = \frac{16}{n_0! \chi} e^{n_0} (2n_0+1+2 \chi) (n_0+2)(4n_0+4)^{n_0+1} \left( en_0/2 \right)^{\chi},
$$
$$
C_0 = \log \left( e^{4.4 n_0+7} n_0^{5.5} D^2 \log ( e D) \right) \; \mbox{ and } \; W_0 = \log \left(
1.5 e B D \log (eD) \right).
$$
Then, if $\log \alpha_1, \ldots, \log \alpha_{n_0}$ are linearly independent over $\mathbb{Z}$ and $b_{n_0} \neq 0$, we have
$$
\log \left| \Lambda \right| > - C(n_0) \, C_0 \, W_0 \, D^2 \, \Omega.
$$
\end{thm}

From (\ref{eqn:main}), we have that
$$
\sqrt{5} y^p - \varepsilon^{2n} = 2 \sqrt{5} - \overline{\varepsilon}^{2n}
$$
and so
\begin{equation} \label{upper-bound}
0 < \Lambda = p \log y + \log ( \sqrt{5}) - 2n \log \varepsilon < \frac{ 2 \sqrt{5}}{ \varepsilon^{2n}} < \frac{2.1}{y^p}.
\end{equation}
We apply Theorem \ref{Matveev} with
$$
D=2, \; \chi = 1, \; n_0=3, \; b_1=1, \; \alpha_1=\sqrt{5}, \; b_2 = -2n, \; \alpha_2 = \varepsilon, \; b_3=p, \; \alpha_3=y,
$$
where, from Lemma \ref{lem:yge19}, we have $y \geq 19$. We may thus take
$$
A_1 =  \log 5, \;  A_2 = \log \varepsilon, \; A_3 = 2 \log y \; \mbox{ and } \;  B = \max \left\{
\frac{n \log \varepsilon}{\log y}, p \right\} =p.
$$
Since
$$
C_0(3) = 2^{18} \cdot 3^2 \cdot 5 \cdot e^4 < 6.45 \times 10^8, \; \; C_0 = \log \left( e^{20.2}  \cdot 3^{5.5} \cdot 4 \log (4e) \right) < 28.5
$$
and 
$$
W_0 =  \log \left(  3 e p \log (2e) \right) < 2.63+\log p
$$
we may  therefore conclude that
$$
\log \Lambda > - 1.139 \cdot 10^{11} \left(  2.63+\log p \right) \log y.
$$
From  (\ref{upper-bound}), we thus have that 
$$
 p \log y <  1.139 \cdot 10^{11} \left(  2.63+\log p \right) \log y + \log (2.1),
$$
and hence
$$
\frac{p}{2.63+\log p} < 1.139 \cdot 10^{11} + \frac{ \log (2.1)}{\left(  2.63+\log p \right) \log y} < 
1.14 \cdot 10^{11}.
$$
We thus have that  $p <3.6 \times 10^{12}$.

Our immediate goal is to sharpen this inequality by proving that $p < 10^{11}$. We will assume for the remainder of this section  that
\begin{equation} \label{first-bound}
10^{11} \leq p <3.6 \times 10^{12}.
\end{equation}

We begin by  appealing to a sharper but less convenient lower bound for linear forms in three complex logarithms, due to Mignotte  (Proposition 5.1 of \cite{Mig2}).

\begin{thm}[Mignotte] \label{miggy}
Consider three non-zero  algebraic numbers $\alpha_1$, $\alpha_2$
and $\alpha_3$, which are either all real and ${}>1,$ or all complex of modulus
one and all ${}\not=1$. Further, assume that the three numbers $\alpha_1, \alpha_2$ and $\alpha_3$ are either all multiplicatively independent, or that two of the number are multiplicatively independent and the third is a root of unity.
We also consider three positive
rational integers $b_1$, $b_2$, $b_3$ with $\gcd(b_1,b_2,b_3)=1$, and the linear form
$$
   \Lambda = b_2\log \alpha_2-b_1 \log \alpha_1-b_3\log \alpha_3 ,
$$
where the logarithms of the $\alpha_i$ are arbitrary determinations of the logarithm,
but which are all real or all purely imaginary.
Suppose further that
$$
   b_2 |\log \alpha_2| =
 b_1\,\vert \log \alpha_1 \vert+  b_3 \,\vert\log \alpha_3\vert \pm  \vert\Lambda\vert 
$$
and put
$$
d_1 = \gcd(b_1,b_2) \; \mbox{ and } \;   d_3 = \gcd(b_3,b_2).
$$
Let $\rho\ge e := \exp(1)$  be a real number. Let $a_1, a_2$ and $a_3$ be real numbers such that
$$
   a_i \ge  \rho \vert \log \alpha_i \vert
   - \log  \vert \alpha_i\vert +2 D \,{\rm h}\kern .5pt(\alpha_i), \qquad
   i \in \{1, 2, 3 \},
$$
where
$\,D=[\mathbb{Q}(\alpha_1,\alpha_2,\alpha_3) : \mathbb{Q}]\bigm/[\mathbb{R}(\alpha_1,\alpha_2,\alpha_3) : \mathbb{R}]$, and assume further that
$$
\Omega := a_1 a_2 a_3 \geq 2.5 \; \mbox{ and } \; a := \min \{ a_1, a_2, a_3 \} \geq 0.62.
$$
Let $m$ and $L$ be positive integers with $m \geq 3$, $L \geq D+4$  and set
$K = [ m \Omega L ].$
Let $\chi$ be fixed with $0 < \chi \leq 2$ and define
$$
c_1 = \max \{ (\chi m L)^{2/3}, \sqrt{2mL/a} \}, \; 
c_2 = \max \{ 2^{1/3} \, (m L)^{2/3}, L \sqrt{m/a} \}, \;
c_3 = (6 m^2)^{1/3} L,
$$
$$
R_i = \left[ c_i a_2 a_3 \right], \; 
S_i = \left[ c_i a_1 a_3 \right] \; \mbox{ and } T_i = \left[ c_i a_1 a_2 \right],
$$
for $i \in \{ 1, 2, 3 \}$, 
and set
$$
R=R_1+R_2+R_3+1, \; S = S_1+S_2+S_3+1 \; \mbox{ and } \; T = T_1+T_2+T_3+1.
$$
Define
$$
c = \max \left\{ \frac{R}{L a_2 a_3}, \frac{S}{L a_1 a_3}, \frac{T}{L a_1 a_2} \right\}.
$$
Finally, assume that the quantity
$$
\begin{array}{c}
\left( \frac{KL}{2} + \frac{L}{4} - 1 - \frac{2K}{3L} \right) \log (\rho) - (D+1) \log L - 3 g L^2 c \, \Omega  \\ \\
  - D (K-1) \log B - 2 \log K + 2 D \log 1.36
 \end{array}
$$
is positive, 
where
$$
   g={1 \over 4}-{K^2L \over 12RST} \; \mbox{ and } \; 
   B = \frac{e^3 c^2 \Omega^2 L^2}{4K^2 d_1 d_3} \left( \frac{b_1}{a_2}+ \frac{b_2}{a_1} \right) 
    \left( \frac{b_3}{a_2}+ \frac{b_2}{a_3} \right).
$$

\noindent {\bf Then either}
\begin{equation} \label{rups}
\log  \Lambda >  - ( KL + \log ( 3 KL)) \log (\rho),
\end{equation}
or the following condition holds :

\smallskip
\noindent {\bf either} there exist non-zero rational integers $r_0$ and $s_0$ such that
\begin{equation} \label{rups2}
   r_0b_2=s_0b_1
\end{equation}
with
\begin{equation} \label{rups3}
   |r_0|
   \le \frac{(R_1+1)(T_1+1)}{M-T_1}
    \;  \mbox{ and } \; 
   |s_0| 
   \le \frac{(S_1+1)(T_1+1)}{M-T_1},
\end{equation}
where
$$
M =  \max\Bigl\{R_1+S_1+1,\,S_1+T_1+1,\,R_1+T_1+1,\,\chi \; \tau_1^{1/2} \Bigr\}, \; \; 
\tau_1 = (R_1+1)(S_1+1)(T_1+1),
$$
{\bf or}
there exist rational integers  $r_1$, $s_1$, $t_1$ and $t_2$, with
$r_1s_1\not=0$, such that
\begin{equation} \label{rups4}
   (t_1b_1+r_1b_3)s_1=r_1b_2t_2, \qquad \gcd(r_1, t_1)=\gcd(s_1,t_2 )=1,
\end{equation}
which also satisfy
$$
    |r_1s_1|
   \le \gcd(r_1,s_1) \cdot
   \frac{(R_1+1)(S_1+1)}{M-\max \{ R_1, S_1 \}},
$$
$$
   |s_1t_1| \le \gcd(r_1,s_1) \cdot
   \frac{(S_1+1)(T_1+1)}{M-\max \{ S_1, T_1 \}} 
  $$
and
$$
   |r_1t_2| % &
   \le \gcd(r_1,s_1) \cdot
  \frac{(R_1+1)(T_1+1)}{M- \max \{ R_1, T_1 \}}.
$$
Moreover, when $t_1=0$ we can take $r_1=1$, and
when $t_2=0$ we can take $s_1=1$.
\end{thm}

\medskip

We apply this result with 
$$
b_2=p, \; \alpha_2 = y, \;b_1=1, \; \alpha_1 = \sqrt{5}, \; b_3= 2n \; \mbox{ and } \; \alpha_3 =\varepsilon,
$$
so that we may take 
$$
D=2, \; d_1=1, \; d_3 \in \{ 1, p \}, \; a_1 = \frac{\rho+3}{2} \log 5, \; a_2= (\rho+3) \log y
$$
and $a_3= (\rho+1) \log ( \varepsilon )$, whence $a=a_3$. 

Notice that, in our situation,  (\ref{rups2}) becomes the equation
$r_0p = s_0$
from which necessarily $|s_0| \geq p > 10^{11}$, whereby (\ref{rups3}) implies that
\begin{equation} \label{frog1}
\frac{(S_1+1)(T_1+1)}{M-T_1} \geq 10^{11}.
\end{equation}
If instead we have (\ref{rups}), then inequality (\ref{upper-bound}) implies that
\begin{equation} \label{frog2}
p \log y < ( KL + \log ( 3 KL)) \log (\rho) + \log (2.1).
\end{equation}

We will choose $L$, $m,$ $\rho$ and $\chi$ to contradict both (\ref{frog1}) and (\ref{frog2}), whereby we necessarily have (\ref{rups4}). 
Specifically, we set
$$
L=485, \; m = 20, \, \rho = 5.7 \; \mbox{ and } \; \chi=2,
$$
so that 
$$
K = \left[ 20 \cdot 485 \cdot  4.35 \log (5) \cdot 6.7 \log ( \varepsilon ) \cdot  8.7 \log y \right],
$$
whereby
$$
1904870 \log y < K  \leq  1904871 \log y.
$$
We have
$$
c_1 < 721.996, \; c_2 < 1207.96, \; c_3 < 6493.5, 
$$
$$
R_1 < 20252 \log y, \; R_2 < 33883 \log y, \; R_3 < 182142 \log y,
$$
$$
S_1 = 16297, \; S_2 = 27266, \; S_3 =146572,
$$
$$
T_1 < 43977 \log y, \; T_2 < 73576 \log y, \; T_3 < 395514 \log y, 
$$
so that
$$
R <236277 \log y +1, \; S = 190136, \; T < 513067 \log y + 1,
$$
and 
$$
c < 17.37,  \; g <0.244 \; \mbox{ and } \; B < 0.3 \, p^2.
$$
We check that
$$
\left( \frac{KL}{2} + \frac{L}{4} \right) \log (\rho) + 4 \log (1.36) > 8.03 \cdot 10^8 \log y,
$$
while, using that $p < 3.6 \cdot 10^{12}$ and $y \geq 19$,
$$
\left(  1 + \frac{2K}{3L} \right) \log (\rho) +3 \log L +3 g L^2 c \, \Omega 
  + 2 (K-1) \log B + 2 \log K < 8.021 \cdot 10^8 \log y.
$$
It follows that the hypotheses of Theorem \ref{miggy} are satisfied. Since we may check that $M > 7.6 \cdot 10^6 \log y$, we have that
$$
\frac{(S_1+1)(T_1+1)}{M-T_1} < 95
$$
contradicting (\ref{frog1}). Also, 
$$
 ( KL + \log ( 3 KL)) \log (\rho) + \log (2.1) < 5 \cdot 10^9,
 $$
contradicting (\ref{frog2}).

We may thus conclude that there exist rational integers  $r_1$, $s_1$, $t_1$ and $t_2$, with
$r_1s_1\not=0$, such that
\begin{equation} \label{pood}
(t_1+2nr_1) s_1=r_1t_2 p,
\end{equation}
where, again using that $M > 7.6 \cdot 10^6 \log y$,
$$
    \left| \frac{r_1s_1}{\gcd(r_1,s_1)} \right|
   \le 43
    \; \mbox{ and } \;
      \left| \frac{s_1t_1}{ \gcd(r_1,s_1)} \right| \le 94.
$$

Since, in all cases, we assume that $p > 10^{11}$, we thus have
$$
\max \{ |r_1|, |s_1|, |t_1| \} < p,
$$
whence, from the fact that  $\gcd(r_1, t_1)=\gcd(s_1,t_2 )=1$, we have $r_1 = \pm s_1$ and so
$t_1+ 2 r_1 n = \pm t_2 p$. Without loss of generality, we may thus write 
$$
u + 2 r |n| = t p,
$$
where $r = |r_1|$ and $t=|t_2|$ are positive integers, $u = \pm t_1$, $r \leq 43$ and $|u| \leq 94$.

We can thus rewrite the linear form
$$
\Lambda = p \log y + \log ( \sqrt{5}) - 2n \log \varepsilon 
$$
as a linear form in two logarithms
\begin{equation} \label{two-logs}
\Lambda =  p \log \left( \frac{y}{\varepsilon^{t/r}}  \right)+ \log \left( 5^{1/2}  \varepsilon^{u/r} \right).
\end{equation}

We are in position to apply the following sharp lower bound for linear forms in two complex logarithms of algebraic numbers, due to Laurent  (Theorem 2 of \cite{Lau}).  

\begin{thm}[Laurent] \label{laurentlemma} 
Let $ \alpha_{1} $ and $ \alpha_{2}$ be multiplicatively independent algebraic numbers, $ h $, $ \rho $ and $ \mu $ be real numbers with $  \rho > 1 $ and $ 1/3 \leq \mu \leq 1$. Set
\begin{center}
$ \begin{array}{ccc}
\sigma=\dfrac{1+2\mu-\mu^{2}}{2}, & \lambda= \sigma \log \rho, & H= \dfrac{h}{\lambda}+ \dfrac{1}{\sigma},
\end{array} $\\
$ \begin{array}{cc}
\omega=2 \left(1+ \sqrt{1+ \dfrac{1}{4H^{2}} } \right), & \theta=\sqrt{1+ \dfrac{1}{4H^{2}} }+ \dfrac{1}{2H}.
\end{array} $
\end{center}
Consider the linear form $ \Lambda=b_{2}\log \alpha_{2}-b_{1}\log \alpha_{1}, $ where $ b_{1} $ and $ b_{2} $ are positive integers.  Put 
$$
 D= \left[  \mathbb{Q}(\alpha_{1}, \alpha_{2} ): \mathbb{Q} \right]/\left[  \mathbb{R}(\alpha_{1}, \alpha_{2} ): \mathbb{R} \right]  
 $$
 and assume that
$$
h \geq \max \left\lbrace D \left( \log \left( \dfrac{b_{1}}{a_{2}}+ \dfrac{b_{2}}{a_{1}} \right) + \log \lambda +1.75 \right)+0.06, \lambda , \dfrac{D \log 2}{2} \right\rbrace ,  
$$

$$
a_{i} \geq \max \left\lbrace 1, \rho \vert \log \alpha_{i} \vert - \log \vert \alpha_{i} \vert + 2Dh(\alpha_{i}) \right\rbrace \ \ \ \ \ (i=1,2), 
$$
and
$$
a_{1}a_{2} \geq \lambda^{2}.
$$

\noindent Then
\begin{equation} \label{laurentall}
 \log \vert \Lambda \vert \geq -C \left( h+ \dfrac{\lambda}{\sigma} \right)^{2} a_{1}a_{2}- \sqrt{\omega \theta} \left(h + \dfrac{\lambda}{\sigma} \right)- \log \left( C' \left(h+ \dfrac{\lambda}{\sigma} \right)^{2} a_{1}a_{2} \right)
 \end{equation}

\noindent with

$$
C=\dfrac{\mu}{\lambda^{3}\sigma} \left( \dfrac{\omega}{6}+ \dfrac{1}{2} \sqrt{\dfrac{\omega^{2}}{9}+ \dfrac{8\lambda \omega^{5/4} \theta^{1/4}}{3 \sqrt{a_{1}a_{2}}H^{1/2} } + \dfrac{4}{3} \left( \dfrac{1}{a_{1}}+ \dfrac{1}{a_{2}} \right) \dfrac{\lambda \omega }{H} } \right)^{2}
$$
and
$$
C'=\sqrt{ \dfrac{C \sigma \omega \theta}{\lambda^{3} \mu} }.
$$
\end{thm}

We apply this result with
$$
b_1=1, \; b_2 = p, \; \alpha_1 = 5^{1/2}  \varepsilon^{u/r}, \; \alpha_2 =  \frac{\varepsilon^{t/r}}{y},
$$
so that $D=2r$,
$$
h(\alpha_1) \leq \frac{\log 5}{2} + \frac{|u|}{2r} \log \varepsilon \; \mbox{ and } \;  h(\alpha_2) \leq \log y + \frac{t}{2r} \log \varepsilon.
$$
We take $\mu=1$ and $\rho = e^4$, so that $\sigma =1$ and $\lambda = 4$.
From (\ref{two-logs}), inequality (\ref{upper-bound}), $|u| \leq 94$, $1 \leq r \leq 43$, and $p > 10^{11}$, we have that
$$
\left| \log y - \frac{t}{r} \log \varepsilon \right| <  \frac{1}{p} \left( \frac{2.1}{y^p} +  \log \left( 5^{1/2}  \varepsilon^{94} \right) \right) < 10^{-9}
$$
and hence may choose 
$$
a_1 = 2562 \; \mbox{ and } \; a_2 = 6 r \log y + 1.
$$
We have
$$
2r \left( \log \left( \dfrac{1}{6r \log y +1}+ \dfrac{p}{2562} \right) + \log 4 +1.75 \right)+0.06 < 2 r \log p
$$
whence 
$$
h = 2 r \log p
$$
is a valid choice for $h$.
A short computation reveals that
$$
C < 0.029, \; \;  C' < 0.044, \; \; 
$$
and hence from (\ref{laurentall}), $y \geq 19$ and $p > 10^{11}$,
$$
\frac{\log \vert \Lambda \vert}{\log y}  > -1784 \left( r \log p+ 2 \right)^{2} r - 1.39 \left( r \log p + 2  \right)- 3 - 0.7 \log \left( r \log p + 2 \right) - \frac{\log r }{\log y}.
$$
From (\ref{upper-bound}), 
$$
\frac{\log \vert \Lambda \vert}{\log y} < \frac{\log 2.1}{\log y} - p < 0.26 -p
$$
whereby  it follows that
$$
p < 1784 \left( r \log p+ 2 \right)^{2} r + 1.39 \left( r \log p + 2  \right) + 3.26 + 0.7 \log \left( r \log p + 2 \right) + \frac{\log r }{\log y}
$$
and so, since $r \leq 43$, we find that $p < 9.1 \cdot 10^{10}$, contradicting (\ref{first-bound}).

%We may thus choose
%$$
%a_1 = 2r \log 5 + 2 |u| \log \varepsilon  + (\rho-1) \left( \frac{\log 5}{2} + \frac{|u| \log \varepsilon}{r} \right) \
%$$
%and 
%$$
%a_2 = 4r \log y + 2t \log \varepsilon + (\rho-1)  \left| \log \left( \frac{\varepsilon^{t/r}}{y} \right) \right|.
%$$

%----------------------------------------------------------------------------
\section{The Method of Kraus}
%----------------------------------------------------------------------------

We let $M_1$ be given by \eqref{eqn:M1}.
The aim of this section is to prove the following proposition, which improves on Lemma~\ref{lem:M1}.
\begin{prop}\label{prop:congmodp}
Let $(n,y,p)$ be a solution to \eqref{eqn:main}
with $p \ge 5$ prime. Then there is an $m_0 \in \{-2,-1\}$ such that 
\begin{enumerate}
\item[(i)] $\overline{\rho}_{E_n,p} \sim \overline{\rho}_{E,p}$
where $E=E_{m_0}$;
\item[(ii)] $n \equiv m_0 \pmod{M_1}$;
\item[(iii)] $n \equiv m_0 \pmod{p}$.
\end{enumerate}
\end{prop}
Observe that this proposition improves over Lemma~\ref{lem:M1}
in two ways. First the elliptic curve $E_{2}$, corresponding
to the \lq pseudo-solution\rq\ $F_4+2=5$, is eliminated.
But also we know that $n \equiv m_0 \pmod{p}$. This will allow
us to rewrite our linear form in three logarithms
as a linear form in two logarithms (Section~\ref{sec:twologs})
and deduce a much sharper bound for the exponent $p$.

In view of Section~\ref{sec:threelogs} we need only prove
Proposition~\ref{prop:congmodp} for prime exponents
$5 \le p < 10^{11}$.
Fix $m_0 \in \{-2,-1,2\}$ and suppose $n \equiv m_0 \pmod{M_1}$.
Write $E=E_{m_0}$. By Lemma~\ref{lem:M1} 
we know that $\overline{\rho}_{E_n,p} \sim \overline{\rho}_{E,p}$.
We shall give a computational criterion, modelled on ideas of Kraus
\cite{Kr} (see also \cite[Lemma 7.4]{BMS}) which allows
us, for each $5 \le p <10^{11}$, to deduce a contradiction when
$m_0=2$, and to deduce $n \equiv m_0 \pmod{p}$ when
$m_0=-2$, $-1$.

Let $k$ be a positive integer satisfying the following:
\begin{enumerate}
\item[(I)] $q=kp+1$ is a prime that is $\equiv \pm 1 \pmod{5}$.
\end{enumerate}

Let $\theta_i$ be the two square roots of $5$ modulo $q$. Then
$q \OO_K=\fq_1 \fq_2$ where the prime ideals $\fq_i$
are given by $\fq_i=(q,\sqrt{5}-\theta_i)$. 
Observe that 
\[
\OO_K/\fq_1=\F_q=\OO_K/\fq_2.
\]
Moreover,
\[
\sqrt{5} \equiv \theta_i \pmod{\fq_i}.
\]
Write
\[
\varepsilon_1=(1+\theta_1)/2.
\]
Then
\[
\varepsilon \equiv \varepsilon_1 \pmod{\fq_1}.
\]
As $\theta_2=-\theta_1$, we know that
\[
\varepsilon \equiv -1/\varepsilon_1 \pmod{\fq_2}.
\]

If $q \mid y$
then $E_n$ has multiplicative reduction at both $\fq_1$
and $\fq_2$. In this case by Lemma~\ref{lem:congs2}
we know that $a_{\fq_i}(E) \equiv \pm (q+1) \equiv \pm 2 \pmod{p}$.
We impose the following condition:
\begin{enumerate}
\item[(II)] $a_{\fq_1}(E) \not\equiv \pm 2 \pmod{p}$
or $a_{\fq_2}(E) \not\equiv \pm 2 \pmod{p}$.
\end{enumerate}
From condition (II) we have $q \nmid y$.
Let $\varrho$ be a primitive root (i.e.\ a cyclic generator) for $\F_q^*$,
and let $\omega=\varrho^p$.
Let
\[
\cY_{q,p}=\{ \omega^{r} \; : \;   0 \le r \le k-1\}.
\]
Thus $y^p \pmod{q} \in \cY_{q,p}$. The set $\cY_{q,p}$
has size $k$. In practice we choose $k$ to be small so
that (I) and (II) are satisfied. This one of the ideas behind
the method of Kraus.

Now fix $\varpi \in \cY_{q,p}$ and suppose $y^p \equiv \varpi \pmod{q}$.
Note that $\sqrt{5} \equiv \theta_1 \pmod{\fq_i}$. By \eqref{eqn:pp2}
we see that
 $\varepsilon^{2n} \pmod{\fq_1}$ is a root 
(in $\F_q$) of the quadratic polynomial
\[
P_\varpi=T^2+(2-\varpi) \cdot \theta_1 \cdot T-1.
\]
We will write
\[
\cT_{q,p}=\{ t \in \F_q \; : \; 
\text{$P_\varpi(t)=0$ for some $\varpi \in \cY_{q,p}$, and $t$ is a square} \}.
\]
Thus $\varepsilon^{2n} \pmod{\fq_1}$ belongs to $\cT_{q,p}$.
The set $\cT_{q,p}$ has at most $2k$ elements. We will reduce
its size using what we know about $n$. Recall that
$n \equiv m_0 \pmod{M_1}$ and therefore $2n \equiv 2m_0 \pmod{2M_1}$.
Let $v=(q-1)/\gcd(q-1,2 M_1)$. It follows that 
$(\varepsilon^{2n}/\varepsilon^{2m_0})^v \equiv 1 \pmod{q}$.
We deduce that $\varepsilon^{2n} \pmod{\fq_1}$ belongs to
\[
\cS_{q,p}(m_0)=\{ t \in \cT_{q,p} \; : \; (t/\varepsilon_1^{2m_0})^v \equiv 1 \pmod{q}\}.
\]

\begin{lem}\label{lem:cRqp}
With notation and assumptions as above, let
$q$ be a prime satisfying conditions (I) and (II). 
Let
\[
\cR_{q,p}(m_0)=\{ t \in \cS_{q,p} : a_q(G_t) \equiv a_{\fq_1}(E) \text{ and } a_q(H_t) \equiv a_{\fq_2}(E) \mod{p} \},
\]
where the elliptic curves $G_t/\F_q$ and $H_t/\F_q$ are given by
\[
G_t \; : \; Y^2=X^3+2(t+\theta_1)X^2+X, \qquad 
H_t \; : \; Y^2=X^3+2(t^{-1}+\theta_2)X^2+X.
\]
Then $\varepsilon^{2n} \pmod{\fq_1}$ belongs to $\cR_{q,p}$.
\end{lem}
\begin{proof}
Let $t \in \cS_{q,p}$ satisfy $t \equiv \varepsilon^{2n} \pmod{\fq_1}$. Then
$G_t$ is the reduction of $E_n$ modulo $\fq_1$,
and $H_t$ is the reduction of $E_n$ modulo $\fq_2$.
In particular, $a_{\fq_1}(E_n)=a_q(G_t)$
and $a_{\fq_2}(E_n)=a_q(H_t)$.
But by Lemma~\ref{lem:congs2} we have $a_{\fq_i}(E_n) \equiv a_{\fq_i}(E)$
for $i=1$, $2$.
It follows that $a_q(G_t) \equiv a_{\fq_1}(E) \pmod{p}$
and $a_q(H_t) \equiv a_{\fq_2}(E) \pmod{p}$. 
\end{proof}

Finally we shall need one more assumption on $q$.
\begin{enumerate}
\item[(III)] $\varepsilon_1^{2k} \not \equiv 1 \pmod{q}$.
\end{enumerate}

\begin{lem}\label{lem:preprop}
Let $(n,y,p)$ be a solution to \eqref{eqn:main} with $p \ge 5$.
\begin{enumerate}
\item[(a)] Let $q$ be a prime satisfying (I), (II) and suppose
that $\cR_{q,p}(2)=\emptyset$. Then $n \not \equiv 2 \pmod{M_1}$.
\item[(b)] Let $q$ be a prime satisfying (I), (II), (III).
Suppose $n \equiv m_0 \pmod{M_1}$ where $m_0=-2$ or $-1$.
Suppose every $t \in \cR_{q,p}(m_0)$ satisfies
$(t/\varepsilon_1^{2m_0}) \equiv 1 \pmod{q}$.
Then $n \equiv m_0 \pmod{p}$.
\end{enumerate}
\end{lem}
\begin{proof}
Part (a) follows immediately from the above. For part (b), 
recall that $\varepsilon \equiv \varepsilon_1 \pmod{\fq_1}$
and also that the reduction of $\varepsilon^{2n}$ modulo $\fq_1$
belongs to $\cR_{q,p}(m_0)$. From the hypothesis
in (b), we have $\varepsilon_1^{2(n-m_0)} \equiv 1 \pmod{q}$.
However $q=2kp+1$ and by assumption (III), $\varepsilon_1^{2k} \not \equiv 1 \pmod{q}$.
Thus $p \mid (n-m_0)$ as required.
\end{proof}

\subsection{Proof of Proposition~\ref{prop:congmodp}}
In Section~\ref{sec:threelogs} we showed that if $n \ne -2$, $-1$, then
$p<10^{11}$. We may therefore assume this bound. We wrote a short \texttt{Magma} script
which for each prime in the range $5 \le p < 10^{11}$, searches for primes $q$
satisfying (I), (II), (III) and applies the criteria in Lemma~\ref{lem:preprop} 
to prove Proposition~\ref{prop:congmodp}.
The total processor time for
the proof is roughly 1200 days, although the computation, running on a 2499MHz
AMD Opterons, was spread over 50 processors, making the actual computation time about 24 days.
%%%%%%%%%%%%%%%%%%%%%%%%%%%%%%%%%%%%%%%%%%%%%%%%%%%%%%%%%%%%%%%%%%%%%%%%%%%%%%%%%%%%%%%%%%%%%

\section{Linear Forms in Two Logs}\label{sec:twologs}

\begin{lem}\label{lem:psmall}
Let $(n,y,p)$ be a solution to \eqref{eqn:main} with $n \ne -1$, $-2$.
%Suppose $n \equiv -1$, $-2$, $2 \pmod{p}$.
Then $p < 5000$.
\end{lem}

Let us assume that $p > 5000$.
Note that $F_{-2n}=-F_{2n}$. Let $N=\lvert n \rvert$ and $Y=\lvert y \rvert$.
Thus
$F_{2N} \pm 2 = Y^p$.
This can be rewritten as
\[
\frac{\varepsilon^{2N}}{\sqrt{5}Y^p}-1=\frac{\varepsilon^{-2N} \mp 2 \sqrt{5}}{\sqrt{5} Y^p}.
\]
Let 
\[
\Delta=2N \log \varepsilon-\log\sqrt{5} -p \log Y.
\]
Using Lemma~B.2 of \cite{Smart}, we have
\[
\lvert \Delta \rvert < \frac{2.1}{Y^p},
\]
and therefore
\begin{equation}\label{eqn:Deltalb}
\log \lvert \Delta \rvert < \log 2.1-p \log(Y).
\end{equation}
By Proposition~\ref{prop:congmodp} we have $n \equiv -1$, $-2 \pmod{p}$.
Thus we can write $N=kp+\delta$ where $\delta=\pm 1$, $\pm 2$. 
Therefore the linear form in three logarithms $\Delta$ may
now be rewritten as a linear form in two logarithms,
$$
\Delta=p \log(\varepsilon^{2k}/Y)-\log(\sqrt{5}/\varepsilon^{2\delta}).
$$
From (\ref{eqn:Deltalb}), $|\delta| \leq 2$, $Y \geq 19$ and $p > 5000$, we have that
\begin{equation} \label{close}
\left| 2k \log \varepsilon - \log Y \right| < \frac{16}{p} < 0.0032.
\end{equation}

We will apply Theorem \ref{laurentlemma}  with
$$
b_1=1, \; \; b_2=p, \; \;   \alpha_1=\sqrt{5}/\varepsilon^{2\delta},
\;  \; \alpha_2=\varepsilon^{2k}/Y, \; \; \mbox{ and }  D=2.
$$

We have that
$$
h (\alpha_1) \leq \frac{1}{2} \log \left( \frac{15}{2} + \frac{7}{2} \sqrt{5} \right) < 1.365
$$
and
$$
h(\alpha_2) \leq \max \{ \log Y, 2k \log \varepsilon \} < 1.01 \log Y,
$$
whereby we can choose, from (\ref{close}), 
$$
a_1 = (\rho-1) \log (\sqrt{5} \, \varepsilon^4) + 5.46
$$
and 
$$
a_2= 0.0032 \, (\rho-1) + 4.04 \log Y.
$$

\begin{lem}
Suppose $n \ne -1$, $-2$. Then $\alpha_1$, $\alpha_2$ are
multiplicatively independent.
\end{lem}
\begin{proof}
If $\alpha_1$, $\alpha_2$ are multiplicatively dependent then
$y=\pm 5^r$ for some $r$. This contradicts Lemma~\ref{lem:powerof5}.
\end{proof}

We now choose $\rho=23$ and check that, in all cases, inequality (\ref{laurentall}) contradicts (\ref{eqn:Deltalb}). This completes the proof of Lemma \ref{lem:psmall}.

%-------------------------------------------
\begin{comment}
\section{An even worse Frey curve}
%-------------------------------------------

We can also think of \eqref{eqn:pp2} as a $(2,3,p)$.
In this case we associate the solution to
\[
E^\prime_n \; : \; Y^2-2\cdot 3^2 \cdot X + 2^2 \cdot 3 \cdot x.
\]
Again the curve is minimal, and with the irreducibility assumption
we obtain
\[
\cN = (2)^7 \cdot (3)^5 \cdot (\sqrt{5}).
\]
\end{comment}

%-------------------------------------------
\section{Deriving the unit equation}
%-------------------------------------------

Our objective is to obtain a bound for $n$ in terms of $p$.
Towards this objective we reduce \eqref{eqn:main}
to a unit equation.
We start with \eqref{eqn:pp2},
where we recall that $x=\varepsilon^{2n}+\sqrt{5}$. Thus
\begin{equation}\label{eqn:s5s6}
(\varepsilon^{2n}+\sqrt{5}+\sqrt{6})(\varepsilon^{2n}+\sqrt{5}-\sqrt{6})
=\sqrt{5}\cdot \varepsilon^{2n}\cdot y^p.
\end{equation}
Let $K=\Q(\sqrt{5})$ and $K^\prime=K(\sqrt{6})$. Write $\OO$
and $\OO^\prime$
for the rings of integers of $K$ and $K^\prime$; 
these both have class number $1$.
As $\gcd(6,y)=1$ (Lemma~\ref{lem:div}), the two factors on the left-hand
side of \eqref{eqn:s5s6} are coprime in $\OO^\prime$. 
The prime ideal $\sqrt{5} \OO$ splits
as a product of two primes in $\OO^\prime$:
\[
\sqrt{5} \OO^\prime= \varphi_1 \OO^\prime \cdot \varphi_2 \OO^\prime
\]
where
\[
\varphi_1 =-2+\sqrt{5}+\left(\frac{1-\sqrt{5}}{2}\right) \sqrt{6} ,
\qquad
\varphi_2=
-2+\sqrt{5}-\left(\frac{1-\sqrt{5}}{2}\right) \sqrt{6} .
%1/2*(-t + 1)*s + t - 2
\]
Let 
\[
\delta=\sqrt{5}+\sqrt{6},\qquad \mu=5+2\sqrt{6}.
\]
Then $\varepsilon$, $\delta$, $\mu$ is a system of fundamental units
for $\OO^\prime$, and the torsion unit group 
is just $\{\pm 1\}$. It follows that
\begin{equation}\label{eqn:descent}
\varepsilon^{2n}+\sqrt{5}+\sqrt{6}= 
\varepsilon^{a} \cdot \delta^{b} \cdot \mu^{c} \cdot \varphi_i
%\left(-2+\sqrt{5}+s \cdot \left(\frac{1-\sqrt{5}}{2}\right) \sqrt{6} \right) 
\cdot \alpha^p
\end{equation}
for some $i \in \{1,2\}$, $a$, $b$, $c \in \Z$ and $\alpha \in \OO^\prime$.
The exponents $a$, $b$, $c$ matter to us only modulo $p$,
as we can absorb any $p$-th power into the term $\alpha^p$.

We now write $M_2=\lcm(M_1,p)$, where $M_1$ is given by
\eqref{eqn:M1}. By Proposition~\ref{prop:congmodp} we know
that $n \equiv -2$ or $-1 \pmod{M_2}$. 
%Thus $a \equiv n \equiv -2$ or $-1 \pmod{p}$.
\begin{lem}\label{lem:binom}
Let $m_0 \in \{-2,-1\}$. 
Let $(n,y,p)$ be a solution to \eqref{eqn:main} with $p\ge 5$,
and $n \equiv m_0 \pmod{M_2}$.
Then
\begin{equation}\label{eqn:binom}
\varepsilon^{2n}+\sqrt{5}+\sqrt{6}=  
\left(\varepsilon^{2m_0}+\sqrt{5} +\sqrt{6} \right) \cdot \alpha^p
\end{equation}
\begin{comment}
If $n \equiv -2 \pmod{M_2}$ then
\[
\varepsilon^{2n}+\sqrt{5}+\sqrt{6}= \left( 
\left(\frac{7-\sqrt{5}}{2}\right) +\sqrt{6} \right) \cdot \alpha^p
\]
for some $\alpha \in \OO_K$. If $n\equiv -1 \pmod{M_2}$ then
\[
\varepsilon^{2n}+\sqrt{5}+\sqrt{6}= \left(
\left(\frac{3+\sqrt{5}}{2}\right) +\sqrt{6} \right) \cdot \alpha^p
\]
\end{comment}
for some $\alpha \in \OO_K$.
\end{lem}
\begin{proof}
Thus the lemma certainly holds if $n=-2$ or $-1$.
We may therefore suppose
$n\ne -2$, $-1$, whence, by Lemma~\ref{lem:psmall},
that $p < 5000$.
We  observe that
\[
%\left(\frac{7-\sqrt{5}}{2}\right) +\sqrt{6}= 
\varepsilon^{-4}+\sqrt{5}+\sqrt{6}=-1\cdot \varepsilon^{-2} \cdot \mu
\cdot \varphi_1,
%\left(-2+\sqrt{5}+\left(\frac{1-\sqrt{5}}{2}\right)\sqrt{6}\right). 
\]
and
\[
%\left(\frac{3+\sqrt{5}}{2}\right) +\sqrt{6}  =
\varepsilon^{-2}+\sqrt{5}+\sqrt{6}=
\varepsilon^{-1} \cdot \delta
\cdot \varphi_2.
%\left(-2+\sqrt{5}-\left(\frac{1-\sqrt{5}}{2}\right)\sqrt{6}\right). 
\]
Thus, if $m_0=-2$, then we want to show,
in \eqref{eqn:descent}, that $i=1$, $(a,b,c) \equiv (-2,0,1) \pmod{p}$,
and if $m_0=-2$, then we want to show that
$i=2$, $(a,b,c) \equiv (-1,1,0) \pmod{p}$.

Observe that 
\[
\Norm_{K^\prime/K}(\varepsilon)=\varepsilon^2,\qquad
\Norm_{K^\prime/K}(\delta)=-1,\qquad
\Norm_{K^\prime/K}(\mu)=1.
\]
Taking norms on both sides of \eqref{eqn:descent} and comparing with
\eqref{eqn:s5s6} we deduce  
$2a \equiv 2n \pmod{p}$ and so $a \equiv n \pmod{p}$.
As $p \mid M_2$ we have derived the required congruences for
$a$. 

Now let $\fq$ be a prime ideal of $K^\prime$ satisfying the following conditions:
\begin{enumerate}
\item[(i)] $\fq$ has degree $1$; we denote $q$ by the rational prime
below $\fq$, and so $\Norm(\fq)=q$.
\item[(ii)] $p \mid (q-1)$.
\item[(iii)] $(q-1) \mid M_2$.
\end{enumerate}
Fix a choice a primitive root $\varpi$ for $\F_\fq^*=\F_q^*$
and we let $\log_{\fq} : \F_\fq^* \rightarrow \Z/p\Z$
be the composition of the discrete logarithm
$\F_\fq^* \rightarrow \Z/(q-1)\Z$ induced by $\varpi$
with the quotient
map $\Z/(q-1)\Z \rightarrow \Z/p\Z$; it is here
that we make use of condition (ii). Now $n \equiv m_0 \pmod{M_2}$
where $m_0=-2$ or $-1$. By assumption (iii) we have
$\varepsilon^{2n} \equiv \varepsilon^{2m_0} \pmod{\fq}$. 
Applying $\log_\fq$ to \eqref{eqn:descent} we obtain 
\[
b \log_\fq(\delta)+c \log_\fq(\mu) \equiv \log_\fq(\varepsilon^{2m_0}+\sqrt{5}+\sqrt{6})-m_0 \log_\fq(\varepsilon)-\log_\fq \varphi_i
%\left(-2+\sqrt{5}+s \cdot \left(\frac{1-\sqrt{5}}{2}\right)\sqrt{6}\right)
\pmod{p}.
\]
It follows that for each choice of $\fq$ satisfying conditions (i), (ii), (iii),
we obtain a linear congruence for $b$, $c$ modulo $p$. 
We wrote a \texttt{Magma} script which did the following.
For each prime $5 \le p < 5000$, each choice $m_0 \in \{-2,-1\}$
and $i \in \{1,2\}$, the script found five prime ideals $\fq$ satisfying (i), (ii), (iii),
and solved the corresponding linear system of congruences for $b$, $c$.
We found that for $m_0=-2$ the system had precisely one
solution when $i=1$, and that solution is $(b,c) \equiv (0,1) \pmod{p}$,
and no solution when $i=2$. Likewise we found that for $m_0=-1$
the system had precisely one solution when $i=2$, namely $(b,c) \equiv (1,0)
\pmod{p}$, and no solution when $i=1$. This completes the proof.
\end{proof}

%To obtain the binomial Thue equation we write $n=kp+m_0$. We let
%$\beta=\varepsilon^{2k}$. Then the conclusion of Lemma~\ref{lem:binom}
%can be written as
%\[
%\varepsilon^{2m_0} \cdot \beta^p - (\varepsilon^{2m_0}+\sqrt{5}+\sqrt{6})\cdot \alpha^p=\sqrt{5}+\sqrt{6}.
%\]

Next we let 
\[
\kappa=(\varepsilon^{2m_0}+\sqrt{5}+\sqrt{6}) \cdot (\sqrt{6}-\sqrt{5}).
\]
Then we can rewrite \eqref{eqn:binom} as
\begin{equation}\label{eqn:preunit}
-(\sqrt{5}+\sqrt{6}) \cdot \varepsilon^{2n}=1-\kappa \cdot \alpha^p.
\end{equation}
The left hand-side is a unit of $K^\prime$. Let 
\[
\tau_j=1-\zeta^j \cdot \sqrt[p]{\kappa} \cdot \alpha, \qquad \zeta=\exp(2 \pi i/p)
\]
for $j=0,\dotsc,p-1$.
It follows that
$\tau_j$ is a unit in the ring of integers
of $K_j=K^\prime(\zeta^j \cdot \sqrt[p]{\kappa})$. 
Let
\[
\nu_0=\zeta^2-\zeta,\qquad
\nu_1=1-\zeta^2, \qquad
\nu_2=\zeta-1.
\]
We obtain the unit equation
\begin{equation}\label{eqn:unit}
\nu_0 \tau_0+\nu_1 \tau_1 + \nu_2 \tau_2=0.
\end{equation}

\section{A bound for $n$}
In this section we derive a bound for the unknown index $n$
in \eqref{eqn:main}. This bound will follow from the bounds
on the heights of the solutions to the unit equation \eqref{eqn:unit}.
to obtain bounds for solutions to unit equations we closely 
follow \cite{integral}. For this we merely need 
some information about the 
number fields containing these solutions.
Recall $K^\prime=\Q(\sqrt{5},\sqrt{6})$
and $K_j=K^\prime(\zeta^j \cdot \sqrt[p]{\kappa})$.
Let $L=K^\prime(\sqrt[p]{\kappa},\zeta)=K_0(\zeta)$. 
\begin{lem}
\[
[K_j : \Q]=4p, \qquad 
[L:\Q]=
\begin{cases}
40 & \text{if $p=5$}\\
4p(p-1) & \text{if $p>5$}.
\end{cases}
\]
Moreover, the signature of $K_j$ is $(4,2p-2)$.
\end{lem}
\begin{proof}
The element $\kappa \in \OO^\prime$ generates
a prime ideal of norm $5$ or $19$ depending on 
whether $m_0=-2$ or $-1$. Thus $[K_j : K^\prime]=p$, and
so $[K_j:\Q]=4p$ by the tower law. To deduce
the signature we observe that for each of the four
embeddings $\sigma: K^\prime \rightarrow \R$,
there is exactly one real choice for the $p$-th
root of $\sigma(\kappa)$, and $(p-1)/2$ complex
conjugate pairs of such choices.

Next we compute $\Q(\zeta) \cap K_j$. Since $\Q(\zeta)$
has degree $p-1$, which is not divisible by $p$, we see
that $\Q(\zeta) \cap K_j=\Q(\zeta) \cap K^\prime$.
If $p>5$ then the intersection is $\Q$, as the intersection
is unramified at all primes.  If $p=5$ then the intersection
in $\Q(\sqrt{5})$. The assertion regarding $[L:\Q]$
follows.
\end{proof}

We shall need a bound for the absolute discriminant of $K_j$;
such a bound is furnished by the following lemma.
\begin{lem}\label{lem:discUB}
Write $\Delta_{K_j}$ for the  
absolute discriminant of $K_j$.
If $m_0=-2$ then $\Delta_{K_j}$ divides 
$2^{6p} \cdot 3^{2p} \cdot 5^{3p-1} \cdot p^{4p}$.
If $m_0=-1$ then $\Delta_{K_j}$ divides 
$2^{6p} \cdot 3^{2p} \cdot 5^{2p} \cdot 19^{p-1} \cdot p^{4p}$.
\end{lem}
\begin{proof}
%We first compute the relative discriminant ideal
%$\Delta_{K_j/K^\prime}$. 
The absolute discriminant of $K^\prime$ is 
$\Delta_{K^\prime}=14400=2^6 \times 3^2 \times 5^2$.
The extension $K_j/K^\prime$
is generated by a root of the polynomial
$x^p-\kappa$, and hence its relative discriminant ideal
divides the discriminant of this polynomial
which is $\pm p^p \kappa^{p-1}$. 
We now apply the following standard formula 
\cite[Theorem 2.5.1]{Cohen} for the absolute
discriminant
\[
\Delta_{K_j} = \pm \Norm_{K^\prime/\Q} (\Delta_{K_j/K^\prime}) \cdot
\Delta_{K^\prime}^{[K_j : K^\prime]}.
\]
The result follows as $\Norm_{K^\prime/\Q}(\kappa)=-5$ or $19$ depending
on whether $m_0=-2$ or $-1$.
\end{proof}

Recall that we are interested in bounding the heights of the
solutions to the unit equation \eqref{eqn:unit} for $m_0=-2$, $-1$
and $5 \le p < 5000$ prime. For each possible choice of $m_0$
and $p$,  
Lemma~\ref{lem:discUB} gives us an upper bound for the 
absolute value of the discriminant
of  $K_j$. Now \cite[Section 5]{integral},
based on a theorem of Landau, gives a computational method 
for deriving an upper bound for the regulators $R_{K_j}$.
As an illustration, we mention that with $p=4999$ (the largest
prime in our range) the bounds we obtain for $R_{K_j}$
are 
\[
R_{K_j} < 2.2 \times 10^{64529},
\qquad
R_{K_j} < 1.4 \times 10^{66241}
\]
for $m_0=-2$, $-1$ respectively.

We now explain how to obtain a bound for $n$. Proposition 8.1 of \cite{integral}
gives positive numbers $A_1$, $A_2$ (depending on the regulators and
unit ranks of the $K_j$) such that 
\begin{equation}\label{eqn:hunit}
h(\nu_2 \tau_2/\nu_0 \tau_0) \le A_2+A_1 \log(H+\max\{h(\nu_j \tau_j) \; : \; j=0,1,2\}),
\end{equation}
where $H$ is an upper bound for the heights $h(\nu_j)$. 
We shall make repeated use of the following properties
for
absolute logarithmic heights (see for example \cite[Lemma 4.1]{integral}): 
\begin{enumerate}
\item[(i)] if $r$ is an integer and $\beta$ is a non-zero algebraic
number then 
$h(\beta^r)=\lvert r\rvert \cdot h(\beta)$;
\item[(ii)] if $\beta_1,\dotsc,\beta_m$ are algebraic numbers then
\[
h(\beta_1  \cdots \cdot \beta_m) \le h(\beta_1)+\cdots+h(\beta_m),
\qquad
h(\beta_1 + \cdots + \beta_m) \le \log{m}+h(\beta_1)+\cdots+h(\beta_m).
\]
\end{enumerate}

As each $\nu_j$ is a sum of two roots of unity, (ii) gives
us $h(\nu_j) \le \log{2}$, so we can take $H=\log{2}$.
By the definition of logarithmic height
\[
h(\tau_2/\tau_0) = h(\nu_2 \tau_2/\nu_0 \tau_0)
\]
since $\nu_2/\nu_0=\zeta^{-1}$ is a root of unity. 
We let $Y=h(\alpha \cdot \sqrt[p]{\kappa})$. Recall
that $\tau_j=1-\zeta^j \cdot \alpha \cdot \sqrt[p]{\kappa}$.
Thus
\[
\alpha \cdot \sqrt[p]{\kappa}=1+\frac{\zeta^2 -1}{\tau_2/\tau_0-\zeta^2}
\]
so
\[
Y \le 3\log{2}+h(\tau_2/\tau_0).
\]
Observe that 
\[
h(\nu_j \tau_j) \le 2 \log{2} +h(\tau_j)=2\log{2}+h(1-\zeta^j \cdot \alpha 
\cdot \sqrt[p]{\kappa})
\le Y+3\log{2}.
\]
From \eqref{eqn:hunit} we deduce
\begin{equation}\label{eqn:hunit2}
Y \le A_2+3\log{2}+A_1 \log(Y+4 \log{2}).
\end{equation}
By Lemma~9.1 of \cite{integral} we have
\[
Y \le 2 A_1 \log{A_1}+2 A_2+10 \log{2}.
\]
Now that we have obtained a bound for $Y=h(\alpha \cdot \sqrt[p]{\kappa})$
we deduce a bound for $n$. From \eqref{eqn:preunit}
\[
\lvert n \rvert \cdot \log(\varepsilon) \le p Y+\log{2}+\frac{1}{2} \log(\sqrt{5}+\sqrt{6}),
\]
which yields a completely explicit bound for $n$. As an illustration
with $p=4999$, we obtain the bounds
\[
\lvert n \rvert \le 2.57 \times 10^{398775},
\qquad
\lvert n \rvert \le 1.01\times 10^{402199},
\]
respectively for $m_0=-2$, $-1$. The corresponding bounds for smaller
values of $p$ are smaller.

\section{Completing the proof of Theorem~\ref{thm-main}}

\begin{lem}
Suppose $p \le 79$. Then $n = -2$, $-1$.
\end{lem}
\begin{proof}
Proposition~\ref{prop:congmodp} tells us that $n \equiv m_0 \pmod{M_1}$ where $M_1$ is given by
\eqref{eqn:M1}, and $m_0=-2$ or $-1$. In fact
\[
M_1 \approx 7.12 \times 10^{4298}.
\]
We computed the bounds for $n$ for all $p < 5000$. We found that for
$p \le 79$ we have
\[
\lvert n \rvert < 1.14\times 10^{4196}, \qquad
\lvert n \rvert < 2.75 \times 10^{4254},
\]
respectively for $m_0=-2$ or $-1$.
Thus for $p \le 79$, $n=-2$ or $-1$.
\end{proof}

The bounds for $n$ we obtain as in the previous section
are larger than $M_1$ for the remaining values $83 \le p < 5000$.
To complete the proof of Theorem~\ref{thm-main}, 
we will show, for $m_0 \in \{-2,-1\}$ and for each of the remaining $p$,
 the
existence of some $M^\prime$ that is much larger than the corresponding
bound for $n$, and such that $n \equiv m_0 \pmod{M^\prime}$.
For this we shall use a very simple sieve. 
Fix a prime $83 \le p < 5000$ and a value $m_0 \in \{-2,-1\}$.
Let $3 \le \ell <10^4$ be a prime, $\ell \ne p$. Suppose we know that
$n \equiv m_0 \pmod{M^\prime}$, where $M^\prime$ is a large
smooth integer, certainly divisible by $M_1$. Let $r=\ord_\ell(M^\prime)$.
We want to show that $n \equiv m_0 \pmod{\ell^{r+1}}$ and so
$n \equiv m_0 \pmod{\ell M^\prime}$. 
We look for primes 
$q \equiv \pm 1 \pmod{5}$
the form $q=kp\ell^{r+1}+1$, such that
$k p \ell^r \mid M^\prime$ (recall that $M_1 \mid M$ is
divisible by all primes $<10^4$ and so certainly divisible by $p$).
Let $\cQ$ be a (small) set of such primes. 
Let 
\[
\cS:=\{ m_0+t \cdot \ell^r \; : \; 0 \le t \le \ell-1\}.
\]
Then $n \equiv m \pmod{\ell^{r+1}}$ for some unique value $m \in \cS$.
We would like to obtain a contradiction for each possible $m \in \cS$
except for $m=m_0$. It would then follow that $n \equiv m_0 \pmod{\ell^{r+1}}$
and so $n \equiv m_0 \pmod{\ell M^\prime}$ as required. 
Fix $m \in \cS$, $m \ne m_0$. 
As $q=kp\ell^{r+1}+1$ and $kp \ell^r \mid M^\prime$,
the assumptions $n \equiv m \pmod{\ell^r}$
and $n \equiv m_0 \pmod{M^\prime}$ force $n \equiv n_q(m) \pmod{q-1}$
for some unique congruence class $n(m,q)$ (which depends on our choices of
$q$ and $m \in \cS$). Now from Binet's formula,
as $q \equiv \pm 1 \pmod{5}$, we have $F_{2n} \equiv F_{2n(m,q)} \pmod{q}$.
Since $F_{2n}+2=y^p$ we have
\begin{equation}\label{eqn:test}
(F_{2n(m,q)}+2)^{k \ell^{r+1}} \equiv 0 \quad \text{or} \quad 1 \pmod{q}.
\end{equation}
If we find some prime $q \in \cQ$ such that \eqref{eqn:test} fails, then
we will have eliminated that particular value of $m$. Once we 
have eliminated all possibilities all $m \in \cS$, $m \ne m_0$,
we will have deduced $n \equiv m_0 \pmod{\ell M^\prime}$ and we can
replace $M^\prime$ by $\ell M^\prime$.

We wrote a simple \texttt{Magma} script which keeps increasing the
exponents of the primes $3 \le \ell<10^4$, $\ell \ne p$ in $M^\prime$
until $M^\prime$ is sufficiently large to deduce that $n=m_0$.
The total processor time for
the proof is roughly 70 days, although the computation, running on a 2499MHz
AMD Opterons, was spread over 50 processors, making the actual computation time about 1.4 days. This completes the proof of Theorem~\ref{thm-main}.

\end{document}